\documentclass{article}
\usepackage{amsfonts}
\usepackage{amssymb}
\usepackage{amsmath}
\usepackage{hyperref}
\usepackage{geometry}
\usepackage{numinsec}
\usepackage{graphics}
\usepackage{epsfig}

\newtheorem{theorem}{Theorem}[section]

\newtheorem{definition}[theorem]{Definition}

\newtheorem{lemma}[theorem]{Lemma}

\newtheorem{proposition}[theorem]{Proposition}

\newenvironment{proof}[1][Proof]{\noindent\textbf{#1.} }{\ \rule{0.5em}{0.5em}}
\input{tcilatex}
\begin{document}

\title{KBSM of the product of a disk with two holes and $S^{1}$}
\author{M. K. Dabkowski \\
Department of Mathematical Sciences\\
University of Texas at Dallas\\
Richardson, TX 75083\\
mdab@utdallas.edu \and M. Mroczkowski \\
Institute of Mathematics\\
University of Gdansk\\
80-952 Gdansk-Oliwa, ul. Wita Stwosza $57$\\
mmroczko@math.univ.gda.pl}
\date{}
\maketitle

\begin{abstract}
\noindent We introduce diagrams and Reidemeister moves for links in $F\times
S^{1}$, where $F$ is an orientable surface. Using these diagrams we compute
(in a new way) the Kauffman Bracket Skein Modules (\emph{KBSM}) for $%
D^{2}\times S^{1}$ and $A\times S^{1}$, where $D^{2}$ is a disk and $A$ is
an annulus. Moreover, we also find the \emph{KBSM} for the $F_{0,3}\times
S^{1},$ where $F_{0,3}$ denotes a disk with two holes, and thus show that
the module is free.\medskip

\bigskip

\end{abstract}

\section{Introduction}

\noindent Skein modules, as invariants of $3$-manifolds, were introduced by
J. Przytycki \cite{P-3} and V. Turaev \cite{T} in 1987 and have become
important algebraic structures for studying $3$-dimensional manifolds $M^{3}$
and knot theory in $M^{3}$. The Kauffman Bracket Skein Module (\emph{KBSM})
is the most extensively studied skein module. We recall its definition here
for the purposes of further sections.

\begin{figure}[h] 
   \centering
   \includegraphics[scale=0.5]{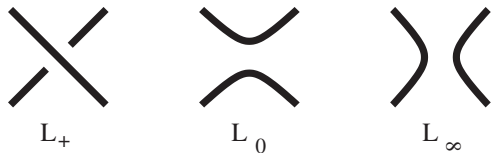} 
   \caption{Smoothings in the Kauffman Bracket skein relation}
   \label{KauffSkeinRel}
\end{figure}

\begin{definition}
\label{def_kbsm}Let $M^{3}$ be an oriented $3$-manifold, $R$ a commutative
ring with identity, and $A\in R$ a unit of $R$. A framed link is an embedded
annulus in which the central curve of the annulus determines an unframed
link. Let $\mathcal{L}_{fr}$ be the set of ambient isotopy classes of
unoriented framed links in $M^{3}$, including the empty link which we denote
by $\emptyset $. We denote by $R\mathcal{L}_{fr}$ the free $R$-module with
basis $\mathcal{L}_{fr}$ $($we choose some ordering of the set $\mathcal{L}%
_{fr})$. Let $S\mathcal{L}_{fr}$ denote the submodule of $R\mathcal{L}_{fr}$
generated by local relations shown 
\begin{align*}
(K1)& : L_{+} =A L_{0} +A^{-1} L_{\infty} \\
(K2)& :\mathcal{L}\sqcup T_{1}=\left( -A^{2}-A^{-2}\right) \mathcal{L}.
\end{align*}%
where $T_{1}$ is the trivial framed link of one component $($the trivial
framed knot$)$ and the triple $L_{+},$ $L_{0},$ $L_{\infty }$ is presented
in \textrm{Figure \ref{KauffSkeinRel}}. Then the \emph{Kauffman Bracket
Skein Module}, $\mathcal{S}_{2,\infty }(M^{3};R,A)$, of $M^{3}$ is defined
to be $\mathcal{S}_{2,\infty }(M^{3};R,A)=R\mathcal{L}_{fr}/S\mathcal{L}%
_{fr}.$
\end{definition}

\noindent The local relations are the ones that arise in the definition of
the Kauffman bracket polynomial of a link. By local relations we mean that
the three links $L_{+}$, $L_{-}$, and $L_{\infty }$ are three framed links
that are identical outside some small neighborhood. It can quickly be seen
that $L^{(1)}=-A^{3}L$ (where $L^{(1)}$ is $L$ with a positive twist) in $%
\mathcal{S}_{2,\infty }(M^{3};R,A) $ which we call the \emph{framing relation%
}.\medskip

\noindent As it was shown in \cite{P-1} and \cite{HP} $\mathcal{S}_{2,\infty
}(S^{3};R,A)$ is free cyclic, $\mathcal{S}_{2,\infty }(S^{1}\times I\times
I;R,A)$ and $\mathcal{S}_{2,\infty }(S^{1}\times S^{1}\times I;R,A)$ are
free modules, generated by infinite sets, whereas $\mathcal{S}_{2,\infty
}(S^{2}\times S^{1};R,A)$ has torsion. However, the \emph{KBSM} has been
found for considerably large classes of $3$\thinspace -\thinspace
dimensional manifolds \cite{P-1, HP, H-P-2, H-P-3, Bu-1}, and the following
problem was proposed in \cite{P} (see Problem $4.4$, pp. 446): \textit{Find 
\emph{KBSM}}\textit{\ of the }$3$\textit{- manifold that is obtained as the
product of a disk with two holes and }$S^{1}$. Let $F_{g,k}$ denote an
orientable surface of genus $g$ with $k$ boundary components. In this paper,
we compute $S_{2,\infty }(F\times S^{1};R,A)$, where $F$ is a disk $%
D^{2}=F_{0,1}$, annulus $A=F_{0,2}$, and disk with two holes $F_{0,3}$. The
modules $S_{2,\infty }(F\times S^{1};R,A)$ for $F=D^{2}$ and $F=F_{0,2} $
were shown to be free in \cite{P-1}, however the methods developed for
computing them will be used in our latter computations of $S_{2,\infty
}(F_{0,3}\times S^{1};R,A).$ Our main result (see Theorem \ref{MainResult})
stating, that $S_{2,\infty }(F_{0,3}\times S^{1};R,A)$ is free solves, in
particular, the Problem $4.4$ of J. Przytycki (\cite{P}, pp. 447) and
supports the conjecture (see Conjecture\textbf{\ }$4.3$, \cite{P}, pp. 446)
that if every closed incompressible surface in $M^{3}$ is parallel to the
boundary of the $3$-manifold $\partial M^{3}$ then $S_{2,\infty }(M^{3};R,A)$
of $M^{3}$ is torsion free\footnote{%
As it was shown in \cite{HP} and \cite{Veve} incompressible (non-boundary
parallel) surface can cause torsion in \emph{KBSM} . Therefore, our result
concerning $S_{2,\infty }(F_{0,3}\times S^{1},%
\mathbb{Z}
\left[ A^{\pm 1}\right] ,A)$ also suggests that torsion in the \emph{KBSM}
cannot be caused by the existence of an immersed torus.}.

\section{Diagrams of links in $F\times S^1$}

\noindent Let $F$ be an orientable surface (possibly with boundary) and $%
I=[0,1]$. Let $f:F\times I\rightarrow F\times S^{1}$ be given by $f(x,y)=(x,$
$e^{2\pi iy})$. Using an argument of general position we can assume that all
links $L$ in $F\times S^{1}$ are transversal to $F\times \{1\}$. Then each $%
f^{-1}(L)$ consists of embedded arcs in $F\times I$ with all endpoints
coming in pairs of the form $\{(x,0),(x,1)\}$ where $x\in F$.

\begin{figure}[h] 
   \centering
   \includegraphics[scale=0.5]{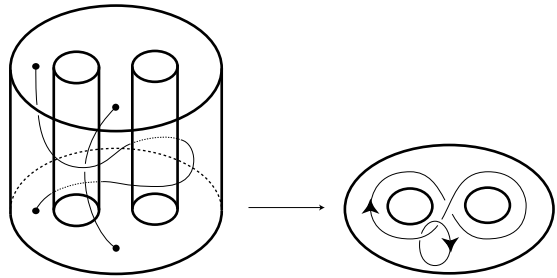} 
   \caption{Links in $F_{0,3}\times S^{1}$}
   \label{diaglink}
\end{figure}


\noindent Let $p$ denote the projection of $F\times I$ onto $F\times \{0\}=F$%
, given by $p(x,y)=(x,0)$. If $D=p(f^{-1}(L))$ then $D$ consists of immersed
curves in $F$. A point in $D$ that corresponds to a pair of endpoints of $%
f^{-1}(L)$ is called a \emph{dot} in $D$. Again, applying an argument of
general position, we may assume that all dots in $D$ are distinct and that
there are only transversal double points in $D$ that are disjoint from dots.
Near each of the double points of $D$ we label two branches as \emph{upper}
and \emph{lower} according to their corresponding values of the arguments $y$
in the second coordinate. Moreover, while the link $L$ crosses the surface $%
f\left( F\times \{1\}\right) =f\left( F\times \{0\}\right) $ (near the point 
$(x_{0},1)\in F\times S^{1}$ of the intersection) the corresponding arc
component containing $(x_{0},1)\in F\times I$ has the $y$ values close to $1$
and increasing in coordinate $y$. Therefore, the corresponding dot in $D$ is
passed in the unique direction (determined by the increasing values of $y)$.
Hence, each link in $F\times S^{1}$ determines uniquely an assignment of
arrows at dots in $D$. Now, a \emph{diagram} of $L$ in $F\times S^{1}$
consists of $D$ together with additional information regarding over and
under branches (for double points of $D$) and an assignment of arrows (for
dots in $D$). An example, showing the construction of the diagram of a link $%
L$ in $F_{0,3}\times S^{1}$ is shown in Figure \ref{diaglink}.\medskip

\begin{figure}[h] 
   \centering
   \includegraphics[scale=0.65]{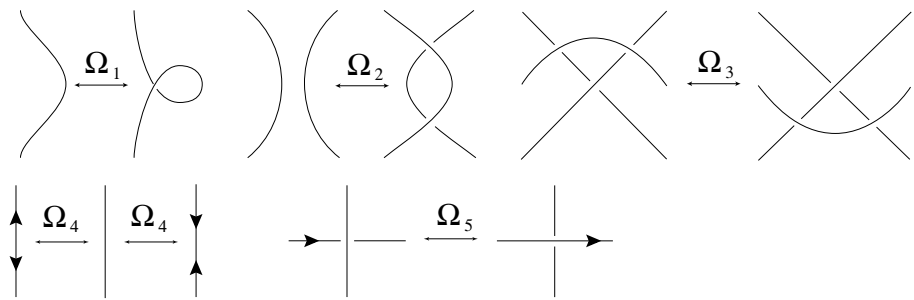} 
   \caption{Reidemeister moves and regular Reidemeister moves}
   \label{reidmv}
\end{figure}

\begin{figure}[h] 
   \centering
   \includegraphics[scale=0.65]{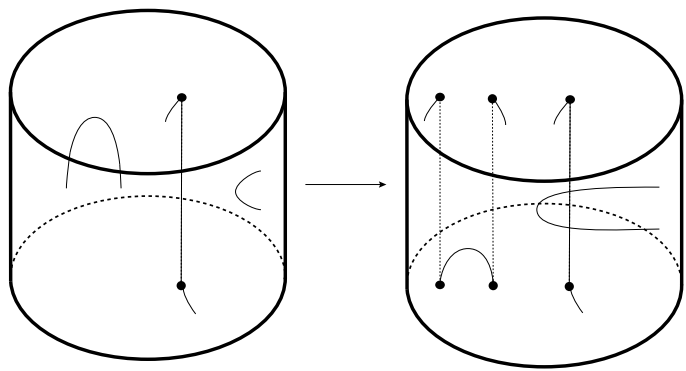} 
   \caption{Geometric interpretation of the $\Omega _{4}$ and $\Omega _{5}$- moves}
   \label{geomint}
\end{figure}

\noindent Two links in $F\times S^{1}$ are isotopic if their diagrams differ
by a finite sequence of "Reidemeister moves". These moves are obtained while
we consider the resolution of generic singularities for the diagrams: cusps,
tangency points, and triple points give us the classical Reidemeister moves $%
\Omega _{1}$, $\Omega _{2}$ and $\Omega _{3}$; double dots give us the
fourth move $\Omega _{4}$; and dots combined with double points give us the
fifth move $\Omega _{5}$. All moves are shown in Figure~\ref{reidmv}. The
geometric interpretation for the $\Omega _{4}$ and $\Omega _{5}$-moves is
shown in Figure~\ref{geomint}. We call the moves $\Omega _{2}$, $\Omega
_{3}, $ $\Omega _{4},$ and $\Omega _{5}$ \emph{regular} Reidemeister moves%
\footnote{%
Analogy with the classical case of links in $S^{3}.$}.

\vspace*{-0.5cm}

\section{Kauffman bracket skein module of $D^{2}\times S^{1}$}

\noindent Let $L$ be a link in $S^{1}\times D^{2}$ and $D$ its diagram. Each
double point of $D$ (crossing of $D$) can be equipped with a positive or a
negative marker corresponding to the horizontal and vertical smoothings.

\begin{figure}[h] 
   \centering
   \includegraphics[scale=0.5]{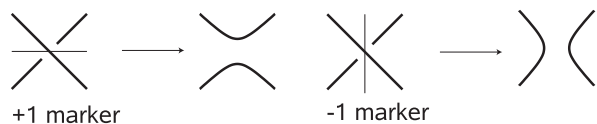} 
   \caption{Positive and negative markers}
   \label{markers}
\end{figure}

\begin{figure}[h] 
   \centering
   \includegraphics[scale=0.5]{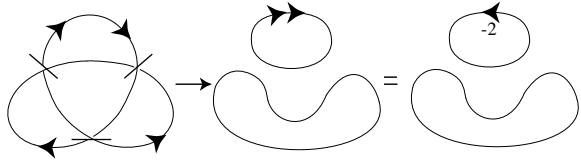} 
   \caption{Diagram with arrows and markers}
   \label{trefoil}
\end{figure}

\noindent A state $s$ of the diagram $D$ is the choice of a marker for each
crossing of $D$. Let $p(s)$ (resp. $n(s)$) be the number of crossings with
positive (resp. negative) markers in $s$. Let $D(s)$ be the diagram obtained
from $D$ by smoothing all crossings in $D$ according to the markers
determined by $s$ and removing all pairs of opposite arrows by the $\Omega
_{4}$-move. To each component of $D(s)$ is assigned an integer that gives
the number of arrows on this component (arrows giving the counterclockwise
orientation of the component are counted as positive, and those giving the
clockwise orientation are counted as negative) as shown in Figure \ref{trefoil}. 
We will refer to such a component as \emph{component with }$n$ \emph{\ arrows}.

\noindent A component $c$ of $D(s)$ is called \emph{trivial} if there are no
arrows on it and every connected component lying inside the disk bounded by $c$ also has no arrows. Let us denote by $|s|$ the number of trivial
components of $D(s)$ and by $D^{\prime }(s)$ the diagram obtained from $D(s)$
by removing all trivial components.

\begin{definition}
\label{def_kb}The Kauffman bracket of $D$ is given by the following sum
taken over all states $s$ of $D:$
\begin{equation*}
\left\langle D\right\rangle
=\sum_{s}A^{p(s)-n(s)}(-A^{2}-A^{-2})^{|s|}D^{\prime }(s)
\end{equation*}
\end{definition}

\begin{lemma}
\label{lemma_234}The Kauffman bracket is preserved by $\Omega _{2}$, $\Omega
_{3},$ and $\Omega _{4}$ moves.
\end{lemma}

\noindent 
\begin{proof}%
The proof for the invariance of the bracket $\left\langle \hspace{0.2cm}%
\right\rangle $ under $\Omega _{2}$ and $\Omega _{3}$-moves is analogous to
the classical case $($Kauffman bracket for classical diagrams \cite{K}$),$
and invariance of $\left\langle \hspace{0.2cm}\right\rangle $ under $\Omega
_{4}$-move follows directly from its definition.%
\end{proof}%
\medskip

\noindent The invariance of $\left\langle \hspace{0.2cm}\right\rangle $
under the $\Omega _{5}$-move requires more detailed analysis. For this
reason, we first introduce a "refined version" of the bracket which is
unchanged under the $\Omega _{5}$-move. Let us denote by $x$ the diagram 
$\epsfig{file=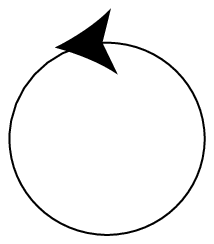, scale=0.2}$ and by $
\mathbb{N}
=\{0,$ $1,$ $2,...\}$ the set of all natural numbers. We show that the $%
S_{2,\infty }(D^{2}\times S^{1};R,A)$ is a free $R$-module, generated by $%
\{x^{n}\mid n\in \mathbb{N}\}$, where $x^{n}$ stands for 
$\epsfig{file=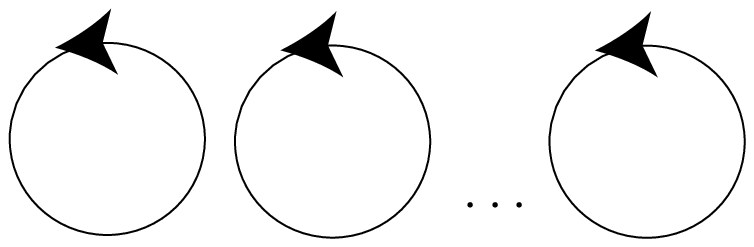, scale=0.2}$ ($n$ copies of $x,$ and $x^{0}=\emptyset $).

\begin{lemma}
\label{lemma_pm}In the skein module $S_{2,\infty }(D^{2}\times S^{1};R,A)$ the following identity holds:
\begin{equation*}
\epsfig{file=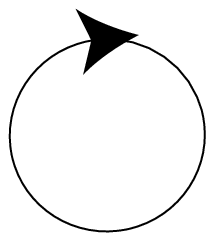, scale=0.2}=A^{-6}\epsfig{file=figx.eps, scale=0.2}
\end{equation*}
\end{lemma}

\noindent 
\begin{proof}%
Using the fact that the removal of a positive kink contributes $-A^{3}$ in
the \emph{KBSM}, and the removal of a negative kink contributes $-A^{-3}$
(just as it is for the classical case of the Kauffman bracket \cite{K}), we
have:%
\begin{equation*}
\epsfig{file=figxinv.eps,scale=0.25}=-A^{-3}\epsfig{file=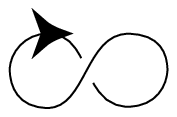,scale=0.55}=-A^{-3}\epsfig{file=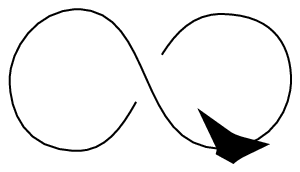,scale=0.3}=A^{-6}\epsfig{file=figx.eps,scale=0.25}.
\end{equation*}
The above calculation finishes our argument.
\end{proof}%
\medskip

\begin{figure}[h] 
   \centering
   \includegraphics[scale=0.6]{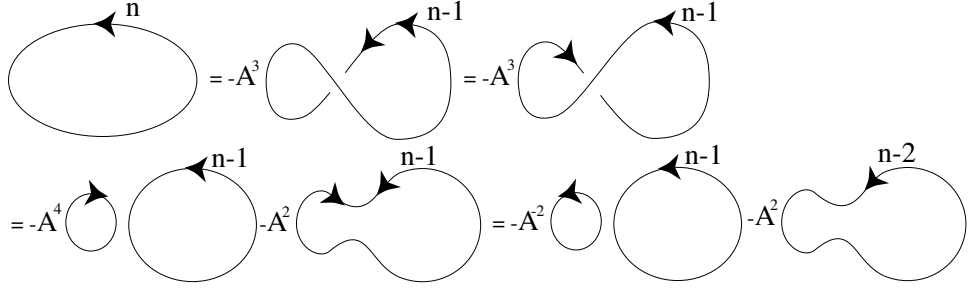} 
   \caption{$D_{n}$ as a combination of $D_{n-1}$ and $D_{n-2}$}
   \label{RecDn}
\end{figure}

\noindent Let $D_{n}$ be the diagram that consists of a single component
with no crossings and $n$ arrows. Using the calculation shown in Figure \ref{RecDn} and applying Lemma \ref{lemma_pm}, we can express $D_{n}$ in 
$S_{2,\infty }(D^{2}\times S^{1};R,A)$ as the combination of $D_{n-1}$ and $D_{n-2}$. Therefore, the following recursion holds true in $S_{2,\infty}(D^{2}\times S^{1};R,A)$:
\begin{equation}
D_{n}=-A^{-2}xD_{n-1}-A^{2}D_{n-2}  \label{eq1}
\end{equation}

\begin{figure}[h] 
   \centering
   \includegraphics[scale=0.6]{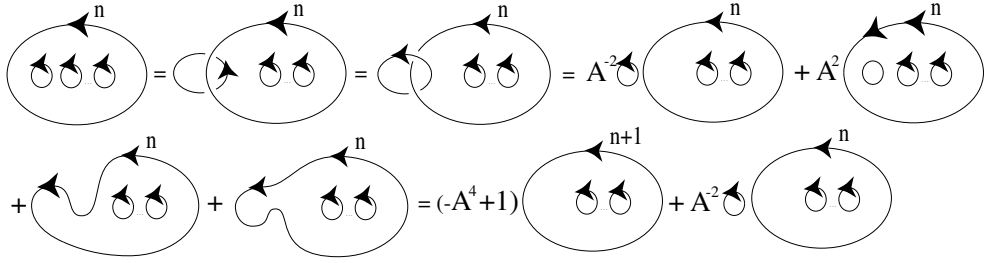} 
   \caption{$D_{n,k}$ as a combination of $D_{n+1,k-1}$ and $D_{n,k-1}$}
   \label{RecDnk}
\end{figure}

\noindent Let $D_{n,k}$ be the
diagram with no crossings and consisting of $k+1$ components, where $k$ of
these components (corresponding to $x^{k})$ are encircled by a single
component with $n$ arrows. We show in Figure \ref{RecDnk} how to express $%
D_{n,k}$ as a combination of $D_{n+1,k-1}$ and $D_{n,k-1}$. Hence, the
following recursion holds true in $S_{2,\infty }(D^{2}\times S^{1};R,A)$:%
\begin{equation}
D_{n,k}=(-A^{4}+1)D_{n+1,k-1}+A^{-2}xD_{n,k-1}  \label{eq2}
\end{equation}

\noindent The recursive relations \ref{eq1} and \ref{eq2} motivate the following definition.

\begin{definition}
Let $P_{n}$ $(n\in \mathbb{Z})$ be polynomials in $x$ and coefficients in
the ring $R$ defined inductively by $P_{0}=-A^{2}-A^{-2},$ $P_{1}=x,$ and
for $n\neq 0,1$ $:$
\begin{equation*}
P_{n}=-A^{-2}xP_{n-1}-A^{2}P_{n-2},
\end{equation*}
where the last relation is also used to define $P_{n}$ for all negative $n$.

\noindent Let $P_{n,k}$ $(n\in \mathbb{Z}$, $k\in \mathbb{N})$ be
polynomials in $x$ and coefficients in $R$ defined inductively by$:$ $P_{n,0}=P_{n},$ and for $k\neq 0:$ 
\begin{equation*}
P_{n,k}=(-A^{4}+1)P_{n+1,k-1}+A^{-2}xP_{n,k-1}
\end{equation*}
\end{definition}

\noindent We have for instance that $P_{-1}=A^{-6}x$. Now we are ready to
define the \emph{refined Kauffman bracket}.

\begin{definition}
\label{def_kbr}Let $D$ be a diagram, $s$ be a state of $D,$ and let $%
D^{\prime }(s)$ denote the corresponding diagram with no crossings and no
trivial components. The \emph{refined Kauffman bracket} polynomial $%
\left\langle D\right\rangle _{r}\in R\left[ x\right] $ is defined
inductively as follows$:$

\begin{description}
\item[$\mathbf{(i)}$] First, we replace by $P_{n}$'s the most nested
components $($no components in the disks they bound$)$ of $D^{\prime}(s)$
that have $n$ arrows, which results in replacing such components by a linear
combination of some $x^{k}$'s.

\item[$\mathbf{(ii)}$] Second, we replace by $P_{n,k}$ each component with $%
n $ arrows that encircles only $x^{k}$.
\end{description}

\noindent The first and second steps are then repeated until $D^{\prime }(s)$
is expressed as a polynomial in $x$ which we denote by $\left\langle
D^{\prime }(s)\right\rangle _{r}$. The \emph{refined} \emph{Kauffman bracket}
of $D$ is given by the following sum taken over all states $s$ of the
diagram $D:$ 
\begin{equation*}
\left\langle D\right\rangle
_{r}=\sum_{s}A^{p(s)-n(s)}(-A^{2}-A^{-2})^{|s|}\left\langle D^{\prime
}(s)\right\rangle _{r}
\end{equation*}
\end{definition}

\noindent We notice that the refined Kauffman bracket was defined in such a
way that it clearly satisfies relations \ref{eq1} and \ref{eq2}. Therefore,
we have the following identities:%
\begin{equation*}
\left\langle D_{n}\right\rangle _{r}=-A^{-2}x\left\langle
D_{n-1}\right\rangle _{r}-A^{2}\left\langle D_{n-2}\right\rangle _{r}\text{
and }\left\langle D_{n,k}\right\rangle _{r}=(-A^{4}+1)\left\langle
D_{n+1,k-1}\right\rangle _{r}+A^{-2}x\left\langle D_{n,k-1}\right\rangle
_{r}.
\end{equation*}
It also follows from the definition that for the refined bracket $\left\langle \hspace{0.2cm}\right\rangle _{r}$, the above two relations hold
true even when they occur in a disk outside of which there are components
which are not involved in the three diagrams appearing in the relations. We
show that $\left\langle \hspace{0.2cm}\right\rangle _{r}$ satisfies a
generalized version of the relation \ref{eq1}.\medskip

\begin{figure}[h] 
   \centering
   \includegraphics[scale=0.6]{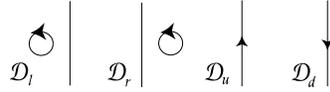} 
   \caption{Diagrams $D_{l}$, $D_{r}$, $D_{u}$ and $D_{d}$}
   \label{figxinout}
\end{figure}

\noindent Let $D_{r}$, $D_{l}$, $D_{u}$, and $D_{d}$ be four diagrams which
are the same outside a small disk but which differ inside the disk as it is
shown in Figure \ref{figxinout}.

\begin{lemma}
\label{lemma_xinout}The refined Kauffman bracket satisfies$:$

\begin{description}
\item[$(1)$] $\left\langle D_{u}\right\rangle _{r}=-A^{-2}\left\langle
D_{r}\right\rangle _{r}-A^{2}\left\langle D_{d}\right\rangle _{r}$

\item[$(2)$] $\left\langle D_{u}\right\rangle _{r}=-A^{-4}\left\langle
D_{l}\right\rangle _{r}-A^{-2}\left\langle D_{d}\right\rangle _{r}$
\end{description}
\end{lemma}

\noindent 
\begin{proof}%
We observe that it suffices to prove $(1)$ and $(2)$ for diagrams with no
crossings and no trivial components. In $D^{\prime }(s)$ let us consider
only the components of $D^{\prime }(s)$ that are inside the component $C$ to
which the vertical segment shown in Figure \ref{figxinout} belongs. Suppose
that $\left\langle \hspace{0.2cm}\right\rangle _{r}$ is applied to these
components so that they are replaced by the linear combination of $x^{k}$'s
(polynomials in $x$). Therefore, it suffices to prove $(1)$ and $(2)$ in the
case when there is only one $x^{k}$ inside $C$.\medskip

\noindent $(1)$ First, let us assume that in $D_{r}$ the circle $x$ is
outside relative to the component $C$. When computing $\left\langle \hspace{0.2cm}\right\rangle _{r}$ any circles inside $C$ are first pushed outside
using the relation \ref{eq2}. This modifies $D_{u}$, $D_{r}$ and $D_{d}$ in
the same way and by the same linear combinations. Thus, we can assume that
there are no $x$-components inside $C$. Considering only $C$ and the circles
appearing in Figure \ref{figxinout} (since contribution of the others is the
same when computing $\left\langle \hspace{0.2cm}\right\rangle _{r}$), we
notice that
\begin{equation*}
D_{u}=D_{n},\text{ }D_{d}=D_{n-2}\text{ and }D_{r}=xD_{n-1},\text{ for some }%
n\in \mathbb{Z}\text{.}
\end{equation*}
Since the refined bracket $\left\langle \hspace{0.2cm}\right\rangle _{r}$
satisfies relation \ref{eq1}, we have:
\begin{equation*}
\left\langle D_{n}\right\rangle _{r}=-A^{-2}x\left\langle
D_{n-1}\right\rangle _{r}-A^{2}\left\langle D_{n-2}\right\rangle _{r},
\end{equation*}
which gives $(1)$.\medskip

\noindent Suppose now that in $D_{r}$ the circle $x$ is inside relative to $C $. As before, when computing 
$\left\langle \hspace{0.2cm}\right\rangle_{r} $, one pushes the $x$'s out of $C$ except for the circle in $D_{r}$. In
this case we observe that:
\begin{equation*}
D_{u}=D_{n},\text{ }D_{d}=D_{n+2},\text{ }D_{r}=D_{n+1,1}\text{ and }D_{l}=xD_{n+1},\text{ for some }n\in \mathbb{Z}.
\end{equation*}
Relation \ref{eq2} for $\left\langle \hspace{0.2cm}\right\rangle _{r}$ gives:
\begin{equation*}
\left\langle D_{r}\right\rangle _{r}=(-A^{4}+1)\left\langle
D_{d}\right\rangle _{r}+A^{-2}\left\langle D_{l}\right\rangle _{r}
\end{equation*}
which is equivalent to:
\begin{equation*}
-A^{-4}\left\langle D_{l}\right\rangle _{r}=(-A^{2}+A^{-2})\left\langle
D_{d}\right\rangle _{r}-A^{-2}\left\langle D_{r}\right\rangle _{r}.
\end{equation*}
Since $\left\langle \hspace{0.2cm}\right\rangle _{r}$ satisfies relation \ref{eq1}, we also have:
\begin{equation*}
\left\langle D_{d}\right\rangle _{r}=-A^{-2}\left\langle D_{l}\right\rangle
_{r}-A^{2}\left\langle D_{u}\right\rangle _{r},
\end{equation*}
which, in turn, is equivalent to:
\begin{equation*}
\left\langle D_{u}\right\rangle _{r}=-A^{-4}\left\langle D_{l}\right\rangle
_{r}-A^{-2}\left\langle D_{d}\right\rangle _{r}.
\end{equation*}
In the last equation, substituting for $-A^{-4}\left\langle
D_{l}\right\rangle _{r}$ from the equation before, we obtain:
\begin{equation*}
\left\langle D_{u}\right\rangle _{r}=(-A^{2}+A^{-2})\left\langle
D_{d}\right\rangle _{r}-A^{-2}\left\langle D_{r}\right\rangle
_{r}-A^{-2}\left\langle D_{d}\right\rangle _{r}=-A^{2}\left\langle
D_{d}\right\rangle _{r}-A^{-2}\left\langle D_{r}\right\rangle _{r},
\end{equation*}
which proves the identity $(1).$\medskip

\noindent $(2)$ From part $(1)$ one has
\begin{equation*}
\left\langle D_{u}\right\rangle _{r}=-A^{-2}\left\langle D_{r}\right\rangle
_{r}-A^{2}\left\langle D_{d}\right\rangle _{r},
\end{equation*}
which is equivalent to:
\begin{equation*}
\left\langle D_{d}\right\rangle _{r}=-A^{-4}\left\langle D_{r}\right\rangle_{r}-A^{-2}\left\langle D_{u}\right\rangle _{r}.
\end{equation*}
Now, we observe that after rotating the diagrams in Figure \ref{figxinout}
by $\pi $, we have $D_{u}$ is switched with $D_{d}$ and $D_{r}$ is switched
with $D_{l}$. Relabelling them accordingly in the last equation gives $(2)$.
\end{proof}%
\medskip

\begin{proposition}
\label{prop_inv} The refined Kauffman bracket $\left\langle \hspace{0.2cm}%
\right\rangle _{r}$ is unchanged under regular Reidemeister moves.
Therefore, $S_{2,\infty }(D^{2}\times S^{1};R,A)$ is a free $R$-module with
basis $\{x^{n}\mid n\in \mathbb{N}\}$.
\end{proposition}

\noindent 
\begin{proof}%
By the result of Lemma \ref{lemma_234}, the Kauffman bracket $\left\langle 
\hspace{0.2cm}\right\rangle $ is preserved by $\Omega _{2},$ $\Omega _{3},$
and $\Omega _{4}$-moves, thus the same holds true for the refined Kauffman
bracket $\left\langle \hspace{0.2cm}\right\rangle _{r}$. It remains to study
the case when diagrams $D$ and $D^{\prime }$ differ by an $\Omega _{5}$
-move. We may assume that there is only one crossing in these diagrams, with
all others being smoothed using the relation $(K1)$. Moreover, we may assume
that inside the disks bounded by the component with the crossing all other
components are already expressed as polynomials in $x$. Moreover, as the
contribution of components outside the component (on which the move is
performed) is the same for both $D$ and $D^{\prime }$, we can additionally
assume that there are no such components. Thus, we need to check the
invariance of $\left\langle \hspace{0.2cm}\right\rangle _{r}$ in the two
cases shown in Figures \ref{figmove5a} and \ref{figmove5b}.

\begin{figure}[h] 
   \centering
   \includegraphics[scale=0.55]{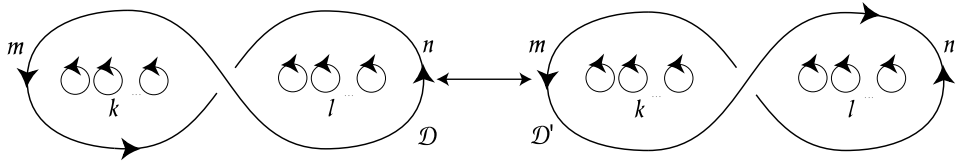} 
   \caption{Invariance of the refined Kauffman bracket}
   \label{figmove5a}
\end{figure}

\noindent In Figure \ref{figmove5a} one needs to verify that $\left\langle
D\right\rangle _{r}=\left\langle D^{\prime }\right\rangle _{r}$. Applying
Lemma \ref{lemma_xinout} for $D$ and $D^{\prime }$, one pushes the $x$'s in
the parts with $m$ arrows (at the expense of getting diagrams with $m+1$ and 
$m-1$ arrows). Analogously, the same happens for the part with $n$ arrows.
Therefore, it is sufficient to check the case when $k=0$ and $l=0$. Now,
applying Lemma \ref{lemma_xinout} again for both diagrams, one may push out
an arrow from the part with $n$ arrows, at the expense of some exterior $x$
and diagrams with $n-1$ and $n-2$ arrows. Thus, it is sufficient to prove
the case for $n=0$ and $n=1$ (one may also reduce $m$ to $0$ or $1$ however,
for our proof it is not needed).\medskip

\noindent For $n=0$, we have:%
\begin{equation*}
\left\langle D\right\rangle _{r}=-A^{-3}\left\langle D_{m+1}\right\rangle
_{r}=-A^{-3}P_{m+1},
\end{equation*}%
and using the inductive definition of $P_{n}$ we have:

\begin{equation*}
\left\langle D^{\prime }\right\rangle
_{r}=AP_{-1}P_{m}+A^{-1}P_{m-1}=A^{-5}xP_{m}+A^{-1}P_{m-1}=-A^{-3}P_{m+1}
\end{equation*}%
For $n=1$, we obtain 
\begin{equation*}
\left\langle D^{\prime }\right\rangle _{r}=-A^{3}\left\langle
D_{m}\right\rangle _{r}=-A^{3}P_{m}
\end{equation*}%
and applying the inductive definition of $P_{n},$ we have:

\begin{equation*}
\left\langle D\right\rangle
_{r}=AP_{m+2}+A^{-1}xP_{m+1}=-A^{-1}xP_{m+1}-A^{3}P_{m}+A^{-1}xP_{m+1}=-A^{3}P_{m}
\end{equation*}%
This finishes our argument for the case shown in Figure \ref{figmove5a}.%

\begin{figure}[h] 
   \centering
   \includegraphics[scale=0.5]{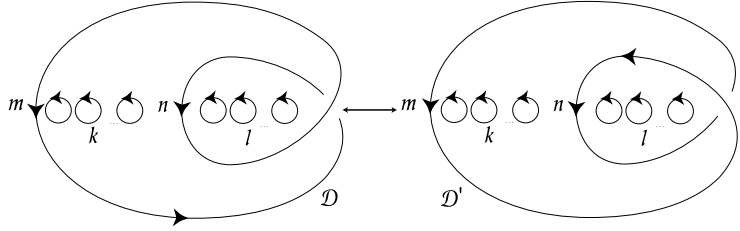} 
   \caption{Invariance of the refined Kauffman bracket}
   \label{figmove5b}
\end{figure}

\noindent In the case shown in Figure \ref{figmove5b}, we also use similar
reasoning as in $(1)$ to reduce it to the case when we can assume that $l=0$
and $k=0$ (applying Lemma \ref{lemma_xinout}). Then again, we apply Lemma %
\ref{lemma_xinout} to push out all the arrows, we arrive again at the
situation when it is sufficient to consider the cases $n=0$ and $n=-1$. We
discuss each one of these cases below.\medskip

\noindent For $n=0$, we have:%
\begin{equation*}
\left\langle D\right\rangle _{r}=-A^{3}P_{m+1}
\end{equation*}%
and using inductive definitions of $P_{n}$ and $P_{n,k}$ we obtain:

\begin{eqnarray*}
\left\langle D^{\prime }\right\rangle _{r}
&=&AP_{m-1}+A^{-1}P_{m,1}=AP_{m-1}+A^{-1}(-A^{4}+1)P_{m+1,0}+A^{-3}xP_{m,0}
\\
&=&AP_{m-1}-A^{3}P_{m+1}-A^{-3}xP_{m}-AP_{m-1}+A^{-3}xP_{m}=-A^{3}P_{m+1}.
\end{eqnarray*}%
For $n=-1$, we have:%
\begin{equation*}
\left\langle D^{\prime }\right\rangle _{r}=-A^{-3}P_{m}
\end{equation*}%
and now applying inductive definitions of $P_{n}$ and $P_{n,k}$ we have:

\begin{eqnarray*}
\left\langle D\right\rangle _{r}
&=&AA^{-6}P_{m+1,1}+A^{-1}P_{m+2}=A^{-5}(-A^{4}+1)P_{m+2,0}+A^{-5}A^{-2}xP_{m+1,0}+A^{-1}P_{m+2}
\\
&=&A^{-5}P_{m+2}+A^{-7}xP_{m+1}=-A^{-3}P_{m}.
\end{eqnarray*}%
The last case verification done for the case shown in Figure \ref{figmove5b}
finishes our argument.%
\end{proof}%
\medskip

\section{Kauffman bracket skein module of $A\times S^{1}$}

\noindent Let $A=F_{0,2}$ denote an annulus. The Kauffman bracket for
diagrams of links in $A\times S^{1}$ is defined in a way that is analogous
to the previous case (see Definition \ref{def_kb}). However, we notice that
after smoothing all of the crossings in such a diagram, there are now two
types of components: ones bounding a disk (called the $x$\emph{-type}) or
ones parallel to the boundary $\partial A=S^{1}\sqcup S^{1}$ components
(called the $y$\emph{-type}).\medskip

\begin{figure}[h] 
   \centering
   \includegraphics[scale=0.5]{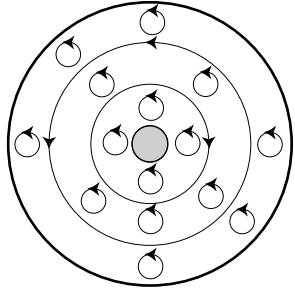} 
   \caption{Diagram corresponding to the word $x^{4}y_{-1}x^{5}y_{2}x^{6}$}
   \label{figxyxy}
\end{figure}

\noindent The Kauffman bracket is first refined in the way analogous to the
one introduced in Definition \ref{def_kbr}. Using the refined bracket for
the $x$-type components enclosed between two successive components of the $y$%
-type allows us to express the $x$-type components as linear combinations of 
$x^{n}$ $(n\geq 0)$ over $R$. Thus, the refined bracket $\left\langle
D\right\rangle _{r}$ for the diagram $D$ can be written as a linear
combination of diagrams shown in Figure \ref{figxyxy}. Moreover, an order of
the two boundary components in $\partial A$ induces a natural order of the $%
y $-type components and all terms $x^{n}$ enclosed by any two successive $y$%
-type components. We use the ordering from the "interior" $S^{1}$ to the
"exterior" $S^{1}$ of $\partial A$ for all of our diagrams. Therefore,
diagrams that constitute terms of the polynomial $\left\langle
D\right\rangle _{r}$ can be encoded uniquely by words of the form

\begin{equation*}
x^{n_{0}}y_{m_{1}}x^{n_{1}}...y_{m_{k}}x^{n_{k}},
\end{equation*}%
where $n_{i}\in \mathbb{N}$, $m_{i}\in \mathbb{Z}$, and $y_{l}$ stands for
the $y$-type component with $l$ arrows (where the sign of $l$ satisfies the
previously introduced convention that counterclockwise orientation is
counted as positive). For example the diagram shown in Figure \ref{figxyxy}
is encoded by the word $x^{4}y_{-1}x^{5}y_{2}x^{6}$. For our convenience, we
identify such diagrams with words.

\noindent Every word $w=y_{m_{1}}y_{m_{2}}...y_{m_{k}}x^{n},$ where $m_{i}=0$
or $1;$ $n\in \mathbb{N}$ will be called \emph{semi-reduced}. For a given
word $w$ we define inductively the semi-reduced refined bracket $%
\left\langle w\right\rangle _{srr}$ corresponding to the diagram of $w$, as
a linear combination of semi-reduced words.\medskip

\noindent $(1)$ If $w$ is semi-reduced, we set $\left\langle w\right\rangle
_{srr}=w$.\medskip

\noindent If $w$ is not semi-reduced then $w$ contains the letter $x$ before
some $y_{n}$ or the letter $y_{n}$ where $n\in 
\mathbb{Z}
\backslash \{0,1\}$. We write $w=w_{1}w_{2},$ where $w_{1}$ is semi-reduced
and contains no $x$.\medskip

\noindent $(2)$ If $w_{2}$ starts with the letter $x$ then $%
w_{2}=x^{k}y_{n}w_{3}$ for some word $w_{3}$. Using the calculation shown in
Figure \ref{RecDnk}, we obtain the following relation in the \emph{KBSM}:%
\begin{equation*}
w_{1}x^{k}y_{n}w_{3}=(-A^{4}+1)w_{1}x^{k-1}y_{n+1}w_{3}+A^{-2}w_{1}x^{k-1}y_{n}xw_{3},
\end{equation*}%
and in this case we define%
\begin{equation*}
\left\langle w\right\rangle _{srr}=(-A^{4}+1)w_{1}\left\langle
x^{k-1}y_{n+1}w_{3}\right\rangle _{srr}+A^{-2}w_{1}\left\langle
x^{k-1}y_{n}xw_{3}\right\rangle _{srr}.
\end{equation*}

\noindent $(3)$ If $w_{2}$ starts with $y_{n}$, where $n>1$ (that is, $%
w_{2}=y_{n}w_{3}$ for some word $w_{3}$)$,$ using the calculation shown in
Figure \ref{RecDn} we obtain the following relation in the \emph{KBSM}: 
\begin{equation*}
w_{1}y_{n}w_{3}=-A^{-2}w_{1}y_{n-1}xw_{3}-A^{2}w_{1}y_{n-2}w_{3},
\end{equation*}%
and accordingly, we define%
\begin{equation*}
\left\langle w\right\rangle _{srr}=-A^{-2}w_{1}\left\langle
y_{n-1}xw_{3}\right\rangle _{srr}-A^{2}w_{1}\left\langle
y_{n-2}w_{3}\right\rangle _{srr}.
\end{equation*}

\noindent $(4)$ If $n<0$ then the last relation is used to express the word
containing $y_{n}$ with words containing $y_{n+1}$ and $y_{n+2}$ and we
define $\left\langle w\right\rangle _{srr}$ accordingly.\medskip

\noindent It can easily be seen that the inductive definition of $%
\left\langle w\right\rangle _{srr}$ results with a linear combination of
semi-reduced words.

\begin{definition}
\label{def_kbsrr}Let $D$ be a diagram and let $\left\langle D\right\rangle
_{r}$ be a linear combination of words $w$ with coefficients in $R$, that
is, $\left\langle D\right\rangle _{r}=\sum_{w}P_{w}w$. We define the \emph{%
semi-reduced refined Kauffman bracket} of the diagram $D$ as follows$:$ 
\begin{equation*}
\left\langle D\right\rangle _{srr}=\sum_{w}P_{w}\left\langle
w\right\rangle_{srr}.
\end{equation*}
\end{definition}

\noindent Now, we prove another generalized version of Lemma \ref%
{lemma_xinout}.

\begin{lemma}
\label{xinout2}The semi-reduced refined Kauffman bracket satisfies the
following properties$:$

\begin{description}
\item[$(1)$] $\left\langle D_{u}\right\rangle _{srr}=-A^{-2}\left\langle
D_{r}\right\rangle _{srr}-A^{2}\left\langle D_{d}\right\rangle _{srr}$

\item[$(2)$] $\left\langle D_{u}\right\rangle _{srr}=-A^{-4}\left\langle
D_{l}\right\rangle _{srr}-A^{-2}\left\langle D_{d}\right\rangle _{srr}$
\end{description}

\noindent where the diagrams $D_{r}$, $D_{l}$, $D_{u},$ and $D_{d}$ are
shown in Figure \ref{figxinout}.
\end{lemma}

\noindent 
\begin{proof}%
If the vertical segment shown in Figure \ref{figxinout} is a part of the $x$%
-type component, then the proof is the same as for the Lemma \ref%
{lemma_xinout} and, moreover, the identities $(1)$ and $(2)$ hold true even
for the refined Kauffman bracket $\left\langle \hspace{0.2cm}\right\rangle
_{r}$. For the other case (when the vertical line segment is a part of the $%
y $-type component) we observe, as before in the proof of Lemma \ref%
{lemma_xinout}, that the relations \ref{eq1} and \ref{eq2} hold true for $%
\left\langle \hspace{0.2cm}\right\rangle _{srr}$, where instead of $D_{n}$
we take $y_{n}$, and any configuration of the remaining components. Indeed,
the components between the first boundary component $S^{1}$ of $\partial A$
and $y_{n}$ are reduced while computing $\left\langle \hspace{0.2cm}%
\right\rangle _{srr}$ to a point when there is some $x^{k}$ to be pushed
through the $y_{n}$. The rest of the proof is the same as for Lemma \ref%
{lemma_xinout}: one pushes all $x$'s, or all but one, through $y_{n}$, and
then $(1)$ follows from definition of $\left\langle \hspace{0.2cm}%
\right\rangle _{srr}$. The property $(2)$ is a consequence of $(1)$ exactly
for the same reasons as in the proof of Lemma \ref{lemma_xinout}.%
\end{proof}%
\medskip

\begin{figure}[h] 
   \centering
   \includegraphics[scale=0.5]{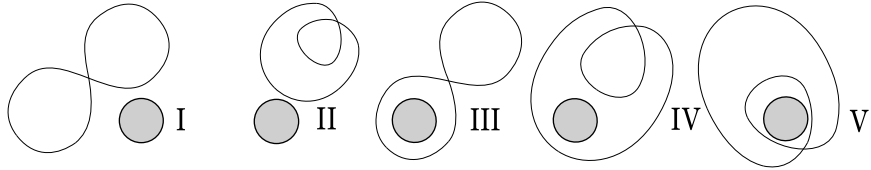} 
   \caption{Five types of $\Omega_5$ moves}
   \label{figfivetypes}
\end{figure}

\noindent Since the Kauffman bracket $\left\langle \hspace{0.2cm}%
\right\rangle $ is invariant under $\Omega _{2}$, $\Omega _{3},$ and $\Omega
_{4}$-moves it follows that the same is true for all other brackets we
defined. To check the invariance of $\left\langle \hspace{0.2cm}%
\right\rangle _{srr}$ under the $\Omega _{5}$-move one may, as before,
consider diagrams where all crossings are smoothed except for the crossing
directly involved in the move. In $S^{1}\times A$ we have five versions
(types) of such moves shown in Figure \ref{figfivetypes}, where only the
first boundary component $S^{1}$ of $\partial A$ is shown and, instead of
presenting the actual $\Omega _{5}$-move, we show only the type of a
component on which this move is performed.

\noindent For the moves of type $\mathbf{I}$ and $\mathbf{II}$ shown in
Figure \ref{figfivetypes} (where smoothing in both ways gives components
bounding disks in $A\,$) $\left\langle \hspace{0.2cm}\right\rangle _{r}$ is
invariant by Proposition \ref{prop_inv} and therefore the same holds true
for $\left\langle \hspace{0.2cm}\right\rangle _{srr}$ and $\left\langle 
\hspace{0.2cm}\right\rangle _{rr},$ where $\left\langle \hspace{0.2cm}%
\right\rangle _{rr}$ is defined after the proof of Proposition \ref%
{prop_invsrr} (see Definition \ref{def_kbrr}).

\begin{proposition}
\label{prop_invsrr} The semi-reduced refined Kauffman bracket, $\left\langle 
\hspace{0.2cm}\right\rangle _{srr}$, is invariant under $\Omega _{5}$-moves
of type $\mathbf{III}$ and $\mathbf{IV}$ shown in \emph{Figure \ref%
{figfivetypes}}.
\end{proposition}

\noindent 
\begin{proof}%
For the move of type $\mathbf{III}$ the situation is shown in Figure \ref%
{figfivea_a}. We may assume that $\left\langle \hspace{0.2cm}\right\rangle
_{srr}$ is applied locally inside of the disks and annuli that are
introduced later after all crossings are smoothed. Applying Lemma \ref%
{xinout2}, we can push all the components $x^{k}$ outside of the component
on which the $\Omega _{5}$-move is applied and then, using the same lemma
again, we can also reduce $n$ to $0$ or $1$. Moreover, we can disregard all
the components inside and outside the component on which $\Omega _{5}$-move
is applied (see Figure~\ref{figfivea_a}) since they equally contribute to $%
\left\langle D\right\rangle _{srr}$ and $\left\langle D^{\prime}\right\rangle _{srr}$.

\begin{figure}[h] 
   \centering
   \includegraphics[scale=0.5]{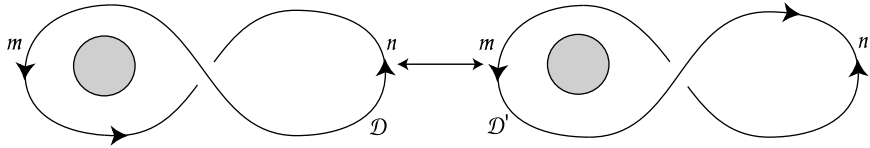} 
   \caption{Invariance of $\left\langle \hspace{0.2cm}\right\rangle _{srr}$ under the $\Omega _{5}$-move of type $\mathbf{III}$}
   \label{figfivea_a}
\end{figure}

\noindent For $n=0$, smoothing both sides gives the following relation%
\begin{equation*}
-A^{-3}y_{m+1}=Ay_{m}A^{-6}x+A^{-1}y_{m-1}
\end{equation*}%
which, after multiplying by $-A^{3}$ gives the relation \ref{eq1} which
holds for $\left\langle \hspace{0.2cm}\right\rangle _{srr}$.\medskip

\noindent For $n=1$, we have:%
\begin{equation*}
Ay_{m+2}+A^{-1}y_{m+1}x=-A^{3}y_{m}
\end{equation*}%
which, after multiplying by $A^{-1}$ becomes again the relation \ref{eq1}
which holds for $\left\langle \hspace{0.2cm}\right\rangle _{srr}$.\medskip

\noindent For type $\mathbf{IV}$ the situation is presented in Figure \ref%
{figfivea_a2}.

\begin{figure}[h] 
   \centering
   \includegraphics[scale=0.5]{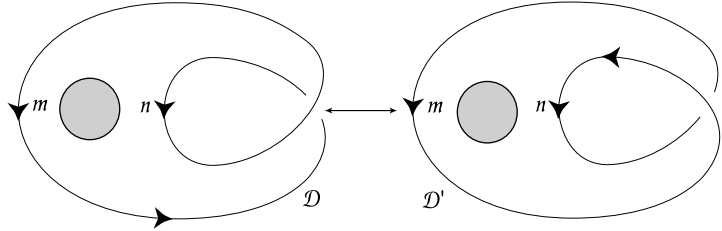} 
   \caption{Invariance of $\left\langle \hspace{0.2cm}\right\rangle _{srr}$ under $\Omega _{5}$-move of type $\mathbf{IV}$}
   \label{figfivea_a2}
\end{figure}

\noindent As in the previous case, we reduce $n$ to $0$ or $-1$ and
disregard all components except the one on which the $\Omega _{5}$-move is
applied. We consider again both cases below.\medskip

\noindent For $n=0$, after smoothing the crossing in $D^{\prime }$ and an
application of relations \ref{eq2} and \ref{eq1}, which hold for $%
\left\langle \hspace{0.2cm}\right\rangle _{srr}$, we obtain%
\begin{eqnarray*}
Ay_{m-1}+A^{-1}xy_{m} &=&Ay_{m-1}+(-A^{3}+A^{-1})y_{m+1}+A^{-3}y_{m}x \\
&=&Ay_{m-1}-A^{3}y_{m+1}-A^{-3}y_{m}x-Ay_{m-1}+A^{-3}y_{m}x=-A^{3}y_{m+1},
\end{eqnarray*}%
which is the same as the result of the smoothing of the crossing in $D$%
.\medskip

\noindent For $n=-1$, smoothing the crossing in $D$ and using relations \ref%
{eq2} and \ref{eq1}, which hold for $\left\langle \hspace{0.2cm}%
\right\rangle _{srr}$, gives%
\begin{eqnarray*}
AA^{-6}xy_{m+1}+A^{-1}y_{m+2}
&=&(-A^{-1}+A^{-5})y_{m+2}+A^{-7}y_{m+1}x+A^{-1}y_{m+2} \\
&=&-A^{-7}y_{m+1}x-A^{-3}y_{m}+A^{-7}y_{m+1}x=-A^{-3}y_{m},
\end{eqnarray*}%
which is the same as the result of smoothing of the crossing in $D^{\prime }$%
.%
\end{proof}%
\medskip

\noindent In order to show the invariance under the $\Omega _{5}$-move of
all five types shown in Figure \ref{figfivetypes}, we need another
refinement where the semi-reduced words are expressed as the \emph{reduced
words} or words of the form $y_{m_{1}}y_{m_{2}}...y_{m_{k+1}}x^{n}$, where
the $m_{1}=m_{2}=...=m_{k}=0$, and $m_{k+1}$ is $0$ or $1$, $n\in \mathbb{N}$%
. To simplify our notations we set $y_{0}=y$ and $y_{1}=y^{\prime },$ thus
the reduced words have the form $y^{k}x^{n}$ or $y^{k}y^{\prime }x^{n}$. Let 
$w$ be a semi-reduced word (and at the same time the diagram represented by
this word). We define inductively the \emph{reduced refinement} $%
\left\langle w\right\rangle _{rr}$ as a linear combination of reduced
words.\medskip

\noindent $(1)$ If $w$ is reduced, we set $\left\langle w\right\rangle
_{rr}=w$.\medskip

\noindent Otherwise $w$ must contain a subword of the form $y^{\prime }y$ or 
$y^{\prime }y^{\prime }$.\medskip

\begin{figure}[h] 
   \centering
   \includegraphics[scale=0.7]{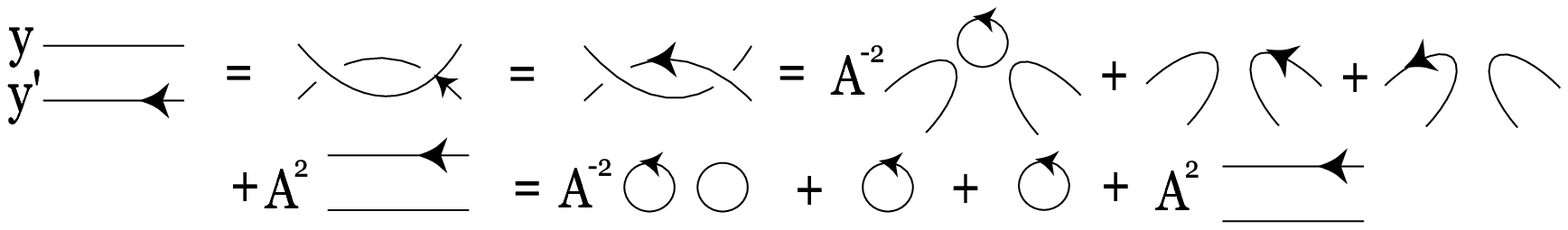} 
   \caption{Reduced refinement $\left\langle w\right\rangle _{rr}$ for $w=y^{k}y^{\prime }yw_{2}$}
   \label{figypy}
\end{figure}

\noindent $(2)$ Suppose that $w=y^{k}y^{\prime }yw_{2}$. In the \emph{KBSM }%
the word $w$ satisfies the following relation:%
\begin{equation*}
y^{k}y^{\prime }yw_{2}=(-A^{-4}+1)y^{k}xw_{2}+A^{2}y^{k}yy^{\prime }w_{2}
\end{equation*}%
as it is shown in Figure \ref{figypy}, and we accordingly define%
\begin{equation*}
\left\langle w\right\rangle _{rr}=(-A^{-4}+1)y^{k}\left\langle \left\langle
xw_{2}\right\rangle _{srr}\right\rangle _{rr}+A^{2}y^{k}y\left\langle
y^{\prime }w_{2}\right\rangle _{rr}
\end{equation*}%
Therefore, the reduced refinement bracket $\left\langle w\right\rangle _{rr}$
can be expressed as the linear combination of the word in which $y$ and $%
y^{\prime }$ are commuted and the word that has one less of the $y$-type
components.

\begin{figure}[h] 
   \centering
   \includegraphics[scale=0.7]{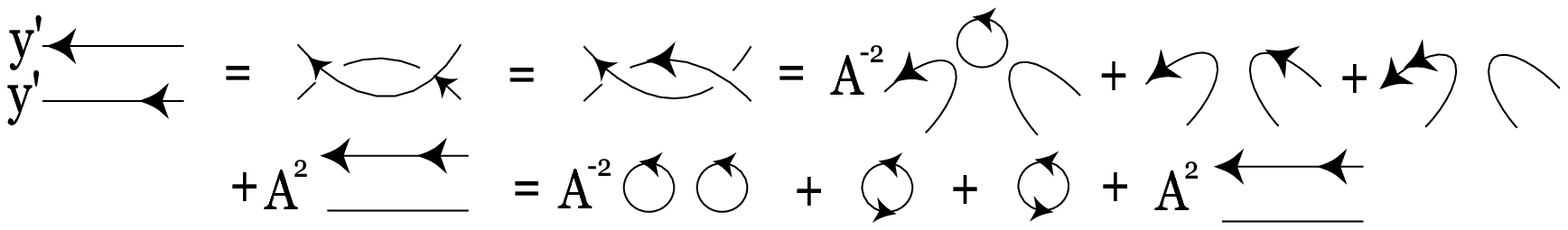} 
   \caption{Reduced refinement $\left\langle w\right\rangle _{rr}$ for $w=y^{k}y^{\prime }y^{\prime }w_{2}$}
   \label{figypyp}
\end{figure}

\noindent $(3)$ Suppose that $w=y^{k}y^{\prime }y^{\prime }w_{2}$. In the 
\emph{KBSM} the word $w$ satisfies the following relation:%
\begin{equation*}
y^{k}y^{\prime }y^{\prime
}w_{2}=A^{-2}y^{k}x^{2}w_{2}+2y^{k}P_{2}w_{2}+A^{2}y^{k}yy_{2}w_{2}
\end{equation*}%
as it is shown in Figure \ref{figypyp}.

\noindent Now we have 
\begin{equation*}
P_{2}=-A^{-2}xP_{1}-A^{2}P_{0}=-A^{-2}x^{2}+A^{4}+1
\end{equation*}%
so $y^{\prime }y^{\prime }=-A^{-2}x^{2}+2A^{4}+2+A^{2}yy_{2},$ and we define

\begin{equation*}
\left\langle w\right\rangle _{rr}=-A^{-2}y^{k}\left\langle \left\langle
x^{2}w_{2}\right\rangle _{srr}\right\rangle
_{rr}+(2A^{4}+2)y^{k}\left\langle w_{2}\right\rangle
_{rr}+A^{2}y^{k}y\left\langle \left\langle y_{2}w_{2}\right\rangle
_{srr}\right\rangle _{rr}
\end{equation*}%
Again, we see that $\left\langle w\right\rangle _{rr}$ can be expressed by
the word containing more $y$'s at the beginning than $w$ (and the number of
components of the $y$-type remains same) and the words containing less $y$%
-type components.\medskip

\noindent It can easily be seen that such an inductive definition of $%
\left\langle w\right\rangle _{rr}$ results in a linear combination of the
reduced words as needed.

\begin{definition}
\label{def_kbrr} Let $D$ be a diagram and let $\left\langle D\right\rangle
_{srr}$ be a linear combination of the semi-reduced words $w$, that is, $%
\left\langle D\right\rangle _{srr}=\sum_{w}P_{w}w$ for some $P_{w}\in R$. We
define the \textit{reduced refined} Kauffman bracket of $D$ by setting%
\begin{equation*}
\left\langle D\right\rangle _{rr}=\sum_{w}P_{w}\left\langle w\right\rangle
_{rr}.
\end{equation*}
\end{definition}

\begin{proposition}
\label{prop_invrr}The reduced refined Kauffman bracket $\left\langle \hspace{0.2cm}\right\rangle _{rr}$ is invariant under all regular Reidemeister
moves. Therefore, $S_{2,\infty}(A\times S^{1};R,A)$ is a free $R$-module
with basis that consists of all reduced words.
\end{proposition}

\noindent 
\begin{proof}%
Being a refinement of the previous bracket, it remains to show that $%
\left\langle \hspace{0.2cm}\right\rangle _{rr}$ is invariant under
Reidemeister $\Omega _{5}$-move of the type $\mathbf{V}$ (see Figure \ref%
{figfivetypes}).The situation is shown in Figure \ref{figfivea_a3}.

\begin{figure}[h] 
   \centering
   \includegraphics[scale=0.5]{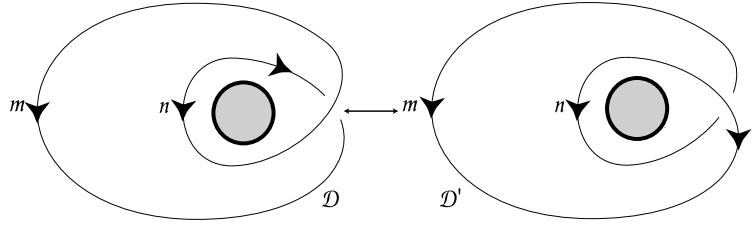} 
   \caption{Invariance of $\left\langle \hspace{0.2cm}\right\rangle _{rr}$ under the Reidemeister moves $\Omega _{5}$ of type $\mathbf{V}$}
   \label{figfivea_a3}
\end{figure}

\noindent Applying Lemma \ref{xinout2} we push out (of the component on
which the $\Omega _{5}$-move is applied) all $x$-type components, and reduce 
$n$ to $0$ or $1$. The components that appear between the component in
Figure \ref{figfivea_a3} and the second boundary component $S^{1}$ of $%
\partial A$ can be disregarded since their contribution to $\left\langle 
\hspace{0.2cm}\right\rangle _{rr}$ is the same for $\left\langle
D\right\rangle _{rr}$ and $\left\langle D^{\prime }\right\rangle _{rr}$ and
comes after the contribution of the component on which the $\Omega _{5}$%
-move is applied. Therefore, we can assume (for the simplicity of our proof)
that there are no such components. For the components of the $y$-type
between the first boundary component $S^{1}$ of $\partial A$ and the
component with the crossing appearing in Figure \ref{figfivea_a3}, we start
applying $\left\langle \hspace{0.2cm}\right\rangle _{rr}$ and arrive at the
situation where these components form $y^{k}$ or $y^{k}y^{\prime }$. Thus,
we have two cases to address:\medskip

\noindent \textbf{Case 1} The $y$-type components form $y^{k}$. Then we have%
\begin{equation*}
\left\langle D\right\rangle _{rr}=A\left\langle
y^{k}y_{n-1}y_{m}\right\rangle _{rr}+A^{-1}\left\langle
y^{k}P_{m-n+1}\right\rangle _{rr}
\end{equation*}%
and 
\begin{equation*}
\left\langle D^{\prime }\right\rangle _{rr}=A\left\langle
y^{k}P_{m-n-1}\right\rangle _{rr}+A^{-1}\left\langle
y^{k}y_{n}y_{m-1}\right\rangle _{rr}.
\end{equation*}%
For simplicity, we omit $y^{k}$ at the beginning of each word. If $n=1$ then
applying Lemma \ref{xinout2} we reduce $m$ to $1$ or $2$.\medskip

\noindent If $n=m=1$ then from the definition of $\left\langle \hspace{0.2cm}%
\right\rangle _{rr}$ it follows that

\begin{equation*}
\left\langle D^{\prime }\right\rangle _{rr}=AP_{-1}+A^{-1}\left\langle
y_{1}y_{0}\right\rangle _{rr}=A^{-5}x+(-A^{-5}+A^{-1})x+Ayy^{\prime
}=Ayy^{\prime}+A^{-1}x=\left\langle D\right\rangle _{rr}
\end{equation*}%
If $n=1$ and $m=2$ then we have:

\begin{eqnarray*}
\left\langle D^{\prime }\right\rangle _{rr}
&=&A(-A^{2}-A^{-2})+A^{-1}\left\langle y^{\prime }y^{\prime }\right\rangle
_{rr}=(-A^{3}-A^{-1})-A^{-3}x^{2}+(2A^{3}+2A^{-1})+A\left\langle
yy_{2}\right\rangle _{rr} \\
&=&A^{3}+A^{-1}-A^{-3}x^{2}+A\left\langle yy_{2}\right\rangle
_{rr}=A\left\langle yy_{2}\right\rangle _{rr}+A^{-1}\left\langle
P_{2}\right\rangle _{rr}=\left\langle D\right\rangle _{rr}
\end{eqnarray*}%
If $n=0$ one reduces $m$ to $0$ or $1$.\medskip

\noindent If $n=m=0$ then using relations \ref{eq1} and \ref{eq2} while
computing $\left\langle \hspace{0.2cm}\right\rangle _{srr}$ we have

\begin{eqnarray*}
\left\langle D\right\rangle _{rr} &=&A\left\langle y_{-1}y\right\rangle
_{rr}+A^{-1}x=-A^{-3}\left\langle yxy\right\rangle _{rr}-A^{-1}\left\langle
y^{\prime }y\right\rangle _{rr}+A^{-1}x \\
&=&-A^{-3}(-A^{4}+1)yy^{\prime }-A^{-5}yyx-A^{-1}(-A^{-4}+1)x-Ayy^{\prime
}+A^{-1}x \\
&=&-A^{-3}yy^{\prime }-A^{-5}yyx+A^{-5}x
\end{eqnarray*}%
and respectively,

\begin{equation*}
\left\langle D^{\prime }\right\rangle _{rr}=A\left\langle
P_{-1}\right\rangle _{rr}+A^{-1}\left\langle yy_{-1}\right\rangle
_{rr}=A^{-5}x-A^{-3}yy^{\prime }-A^{-5}yyx.
\end{equation*}%
Therefore, we have%
\begin{equation*}
\left\langle D\right\rangle _{rr}=\left\langle D^{\prime }\right\rangle
_{rr}.
\end{equation*}%
If $n=0$ and $m=1$ then, we have

\begin{eqnarray*}
\left\langle D\right\rangle _{rr} &=&A\left\langle y_{-1}y_{1}\right\rangle
_{rr}+A^{-1}\left\langle P_{2}\right\rangle _{rr}=-A^{-1}\left\langle
y^{\prime }y^{\prime }\right\rangle _{rr}-A^{-3}\left\langle yxy^{\prime
}\right\rangle _{rr}-A^{-3}x^{2}+A^{3}+A^{-1} \\
&=&A^{-3}x^{2}-2A^{3}-2A^{-1}-A\left\langle yy_{2}\right\rangle
_{rr}-A^{-3}(-A^{4}+1)\left\langle yy_{2}\right\rangle
_{rr}-A^{-5}yy^{\prime }x-A^{-3}x^{2}+A^{3}+A^{-1} \\
&=&-A^{3}-A^{-1}-A^{-3}\left\langle yy_{2}\right\rangle
_{rr}-A^{-5}yy^{\prime }x=-A^{3}-A^{-1}+A^{-5}yy^{\prime
}x+A^{-1}yy-A^{-5}yy^{\prime }x \\
&=&-A^{3}-A^{-1}+A^{-1}yy
\end{eqnarray*}%
and respectively, for $\left\langle D^{\prime }\right\rangle _{rr}$ we have

\begin{equation*}
\left\langle D^{\prime }\right\rangle
_{rr}=A(-A^{2}-A^{-2})+A^{-1}yy=-A^{3}-A^{-1}+A^{-1}yy
\end{equation*}%
Therefore, again it follows that%
\begin{equation*}
\left\langle D\right\rangle _{rr}=\left\langle D^{\prime }\right\rangle _{rr}
\end{equation*}

\noindent \textbf{Case 2} The $y$-type components form $y^{k}y^{\prime }$.
The situation is shown in Figure \ref{figfivea_a4}.

\begin{figure}[h] 
   \centering
   \includegraphics[scale=0.4]{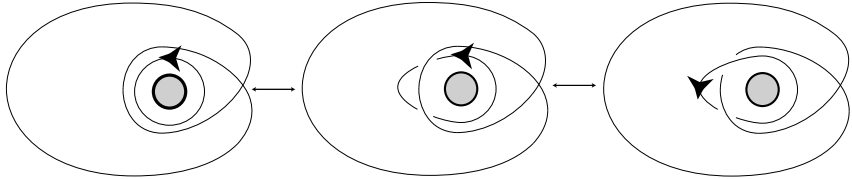} 
   \caption{Invariance of $\left\langle \hspace{0.2cm}\right\rangle _{rr}$ under the Reidemeister moves $\Omega _{5}$ of type $\mathbf{V}$}
   \label{figfivea_a4}
\end{figure}

\noindent We apply the $\Omega _{2}$-move followed by the $\Omega _{5}$-move
on both $D$ and $D^{\prime }$ and observe that the refined bracket $%
\left\langle \hspace{0.2cm}\right\rangle _{rr}$ is unchanged for both
diagrams. This is clear for the $\Omega _{2}$-move. For the $\Omega _{5}$%
-move shown in Figure \ref{figfivea_a4} one verifies the invariance of $%
\left\langle \hspace{0.2cm}\right\rangle _{rr}$ by considering all
smoothings of the two crossings that are not involved in the move. After
doing so, one obtains moves that are not of the type $\mathbf{V}$ and a move
of type $\mathbf{V}$ to which \textbf{Case 1} applies (because inside the
component $y^{\prime }$ there is only $y^{k}$). Now, we apply the desired $%
\Omega _{5}$-move between $D$ and $D^{\prime }$ that we just transformed by
the two Reidemeister moves in Figure \ref{figfivea_a4}. The $\Omega _{5}$%
-move leaves $\left\langle \hspace{0.2cm}\right\rangle _{rr}$ unchanged
because again, by considering all smoothings of the two crossings that are
not involved in this move, we obtain moves that are not of type $\mathbf{V}$
and a move of type $\mathbf{V}$ to which \textbf{Case 1 }applies.%
\end{proof}%

\section{Kauffman bracket skein module of $F_{0,3}\times S^{1}$}

\begin{figure}[h] 
   \centering
   \includegraphics[scale=0.6]{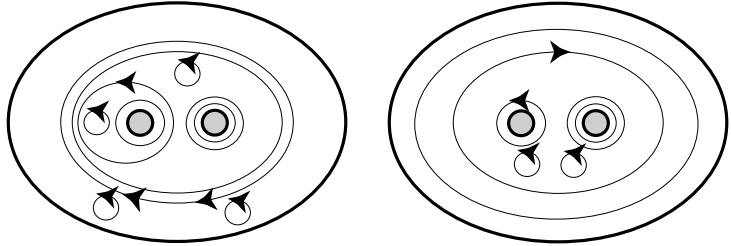} 
   \caption{Diagrams corresponding to the words $y_{0}xy_{1}z_{0}z_{0}x^{2}t_{2}t_{0}x$ and $y^{\prime}z^2tt^{\prime}x^2$}
   \label{figxyzt}
\end{figure}

\noindent Recall that we denoted the disk with two holes by $F_{0,3}$. The
Kauffman bracket, for the case $F_{0,3}\times S^{1},$ is defined again just
as in the Definition~\ref{def_kb}. After smoothing all crossings in a
diagram of the link $L$ we can now encounter four types of components: ones
bounding a disk (of the $x$\emph{-type}) or ones parallel to one of the
three boundary components $S^{1}$ of $\partial F_{0,3}$. We represent $%
F_{0,3}$ as the disk $D^{2}$ with two smaller disks removed, one on the left
denoted by $D_{l}^{2}$ and one on the right denoted by $D_{r}^{2}$.
Components parallel to $\partial D_{l}^{2}$, $\partial D_{r}^{2}$ and $%
\partial D^{2}$ are called, respectively, of the \emph{y-type}, \emph{z-type}
and \emph{t-type}. As in the case of $A\times S^{1}$, diagrams can be
expressed using words. An example of a diagram corresponding to the word $%
y_{0}xy_{1}z_{0}z_{0}x^{2}t_{2}t_{0}x$ is shown on the left in Figure~\ref%
{figxyzt}. Note that for the components of $t$-type the order is from the
circle $\partial D^{2}$ to the interior and the arrows corresponding to the 
\emph{clockwise} orientation are counted as positive. In such words
components of the $y$-type are always written before the components of the $%
z $-type and the $t$-type, and components of the $z$-type are always written
before components of the $t$-type. If $x^{d}$ appears at the end of the
word, it is in between the $y$, $z$ and $t$ components. Otherwise, if it is
before--say the $y$-type component--then it is placed somewhere in between
the $y$ components as it is shown on the left in Figure~\ref{figxyzt}. As
before, to simplify our notations we set $y_{0}=y$, $y_{1}=y^{\prime }$, $%
z_{0}=z$, $z_{1}=z^{\prime }$, $t_{0}=t$ and $t_{1}=t^{\prime }$.

\noindent The bracket $\left\langle \hspace{0.2cm}\right\rangle _{r}$ is
constructed exactly like for $D^{2}\times S^{1}$ and $A\times S^{1}$, and
the constructions of $\left\langle \hspace{0.2cm}\right\rangle _{srr}$ and $%
\left\langle \hspace{0.2cm}\right\rangle _{rr}$ are done analogously as for $%
A\times S^{1}$ by considering separately the components of types $y$, $z$,
and $t$. The $x$-type components and the arrows are pushed by this procedure
towards the area between the components of the three types. In this way $%
\left\langle \hspace{0.2cm}\right\rangle _{r}$ is constructed and allows us
to expresses diagrams as the linear combination of the words $y^{a}y^{\prime
a^{\prime }}z^{b}z^{\prime b^{\prime }}t^{c}t^{\prime c^{\prime }}x^{d}$
where $a,b,c,d\in \mathbb{N}$ and $a^{\prime },b^{\prime },c^{\prime }\in
\{0,1\}$. Such words, as before, are called \emph{reduced}. An example,
represented by $y^{\prime }z^{2}tt^{\prime }x^{2}$, appears on the right in
Figure~\ref{figxyzt}. \medskip

\begin{figure}[h] 
   \centering
   \includegraphics[scale=0.6]{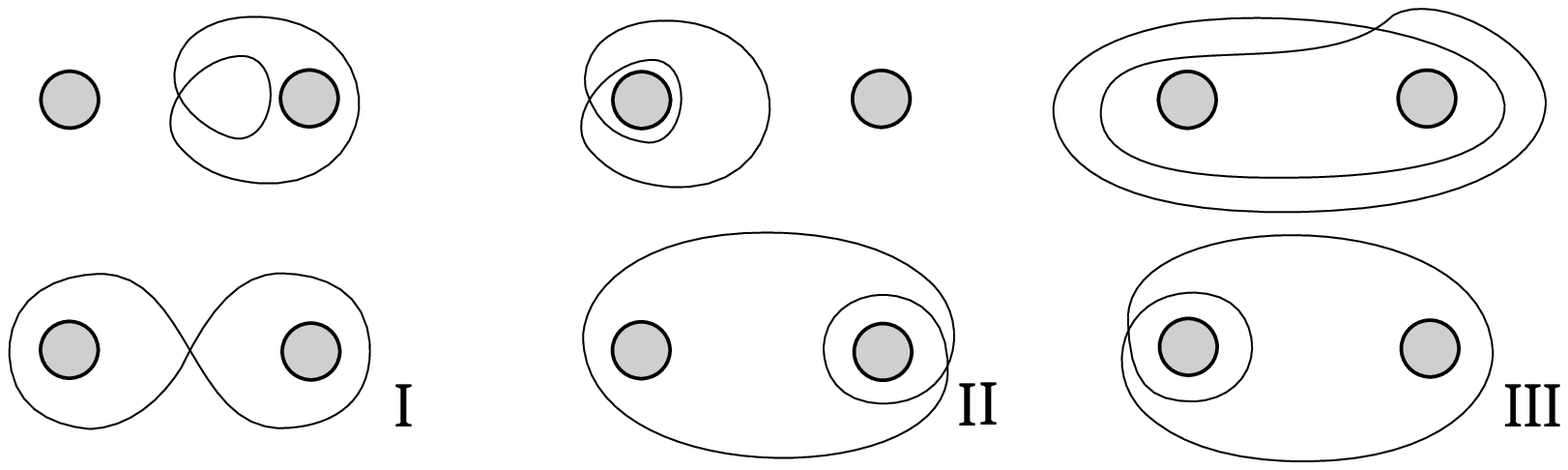} 
   \caption{Moves of types $\mathbf{I},$ $\mathbf{II},$ and $\mathbf{III}$}
   \label{figthreetypes}
\end{figure}

\noindent The refined reduced bracket $\left\langle \hspace{0.2cm}%
\right\rangle _{rr}$ is invariant under $\Omega _{5}$-moves shown in Figure~%
\ref{figfivetypes} (where the additional hole $D^{2}$ is not shown and is
outside of the moves) using the same arguments as in the proofs for $%
D^{2}\times S^{1}$ and $A\times S^{1}$. Figure~\ref{figthreetypes} shows
some of the types of $\Omega _{5}$-moves for which $\left\langle \hspace{%
0.2cm}\right\rangle _{rr}$ is invariant, and three new types of moves, for
which we introduce new refinements of the bracket necessary for us to show
the invariance under the $\Omega _{5}$-move. We first construct the
refinement which will be invariant under the $\Omega _{5}$-move for all
old-type diagrams and some new types of moves. Then, finally, we construct a
refinement which will be invariant under all $\Omega _{5}$-moves, therefore
letting us define the map from $S_{2,\infty }(F_{0,3}\times S^{1};R,A)$ to a
free $R$-module with an explicit basis.

\noindent We notice that in the \emph{KBSM} it is possible for arrows to
jump between components of type $y$, $z$ and $t$. Using this idea, we define
the \emph{quasi-final bracket} $\left\langle w\right\rangle _{qf}$ for the
reduced words $w$. A reduced word is called \emph{quasi-final} if it has at
most one prime (i.e. $y^{\prime }$, $z^{\prime }$ or $t^{\prime }$) and
after the occurrence of such a prime only $x^{k}$ can follow it. For
instance, the reduced word $yzz^{\prime }tx^{2}$ is not quasi-final because
after $z^{\prime }$ we have $t$ that follows it, whereas the reduced word $%
yy^{\prime }x$ is quasi-final. \medskip \noindent Letting $w$ be a reduced
word, we define the quasi-final bracket $\left\langle w\right\rangle _{qf}$
inductively:

\begin{description}
\item[$(1)$] If $w$ is a quasi-final word then we set, $\left\langle
w\right\rangle _{qf}=w$.
\end{description}

\begin{figure}[h] 
   \centering
   \includegraphics[scale=0.65]{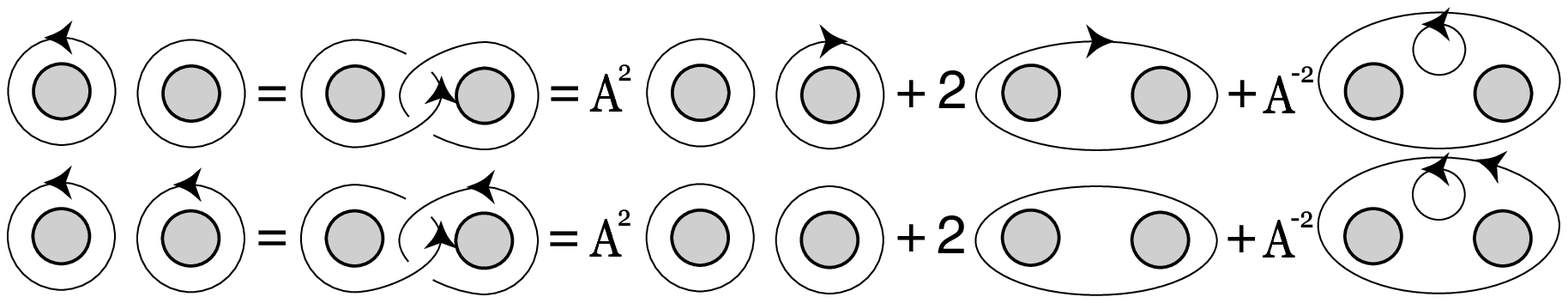} 
   \caption{Quasi-final bracket $\left\langle w\right\rangle _{qf}$ for $w=y^{a}y^{\prime}z^{b}t^{c}t^{\prime c^{\prime }}x^{d}$ (top) and $w=y^{a}y^{\prime}z^{b}z^{\prime }t^{c}t^{\prime c^{\prime }}x^{d}$ (bottom)}
   \label{figqf_yz}
\end{figure}

Otherwise there are several cases to consider depending on the form of the
reduced word $w$ and values for $a,b,c,d\in \mathbb{N}$:

\begin{description}
\item[$(2)$] If $w=y^{a}y^{\prime }z^{b}t^{c}t^{{\prime }c^{\prime }}x^{d}$,
where $c^{\prime }=0$ or $1$ and $b>0,$ then in the \emph{KBSM}, as shown on
the top part of Figure~\ref{figqf_yz}, the word $w$ satisfies the following
relation:

\begin{equation*}
y^{a}y^{\prime }z^{b}t^{c}t^{\prime c^{\prime
}}x^{d}=A^{2}y^{a+1}z^{b-1}z_{-1}t^{c}t^{\prime c^{\prime
}}x^{d}+2y^{a}z^{b-1}t^{c}t^{\prime c^{\prime }}x^{d}t^{\prime
}+A^{-2}y^{a}z^{b-1}t^{c}t^{\prime c^{\prime }}x^{d}tx.
\end{equation*}

\item[$(3)$] If $w=y^{a}y^{\prime }z^{b}z^{\prime }t^{c}t^{{\prime }%
c^{\prime }}x^{d}$, where $c^{\prime }=0$ or $1$, then in the \emph{KBSM},
as shown on the bottom part of Figure~\ref{figqf_yz}, the word $w$ satisfies
the following relation:

\begin{equation*}
y^{a}y^{\prime }z^{b}z^{\prime }t^{c}t^{{\prime }c^{\prime
}}x^{d}=A^{2}y^{a+1}z^{b+1}t^{c}t^{{\prime }c^{\prime
}}x^{d}+2y^{a}z^{b}t^{c}t^{\prime c^{\prime
}}x^{d}t+A^{-2}y^{a}z^{b}t^{c}t^{\prime c^{\prime }}x^{d}t_{-1}x.
\end{equation*}
\end{description}

\begin{figure}[h] 
   \centering
   \includegraphics[scale=0.65]{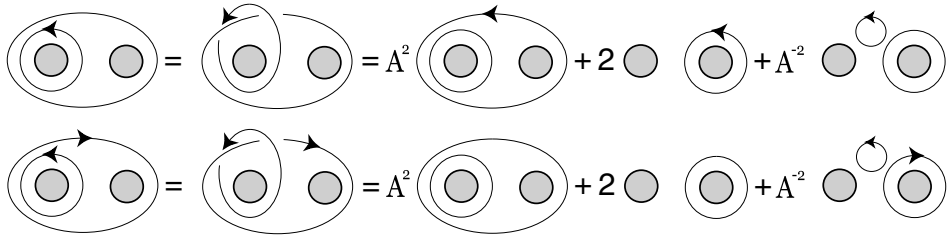} 
   \caption{Quasi-final bracket $\left\langle w\right\rangle _{qf}$ for $w=y^{a}y^{\prime }t^{c}x^{d}$ (top) and $w=y^{a}y^{\prime }t^{c}t^{\prime }x^{d}$ (bottom)}
   \label{figqf_yt}
\end{figure}

\begin{description}
\item[$(4)$] If $w=y^{a}y^{\prime }t^{c}x^{d}$, where $c>0$, then in the 
\emph{KBSM}, as shown on the top of Figure~\ref{figqf_yt}, the word $w$
satisfies the following relation: 
\begin{equation*}
y^{a}y^{\prime
}t^{c}x^{d}=A^{2}y^{a+1}t^{c-1}t_{-1}x^{d}+2y^{a}x^{d}z^{\prime
}t^{c-1}+A^{-2}y^{a}x^{d}zt^{c-1}x.
\end{equation*}

\item[$(5)$] If $w=y^{a}y^{\prime }t^{c}t^{\prime }x^{d}$, then in the \emph{%
KBSM}, as shown on the bottom of Figure~\ref{figqf_yt}, the word $w$
satisfies the following relation: 
\begin{equation*}
y^{a}y^{\prime }t^{c}t^{\prime
}x^{d}=A^{2}y^{a+1}t^{c+1}x^{d}+2y^{a}x^{d}zt^{c-1}+A^{-2}y^{a}x^{d}z_{-1}t^{c-1}x.
\end{equation*}%
\end{description}

\noindent For cases $(6)$ and $(7)$ similar identities to the ones shown in
Figure~\ref{figqf_yt} hold, with the roles of $y$ and $z$ components
switched.

\begin{description}
\item[$(6)$] If $w=y^{a}z^{b}z^{\prime }t^{c}x^{d}$, where $c>0$, then in
the \emph{KBSM } the word $w$ satisfies the following relation: 
\begin{equation*}
y^{a}z^{b}z^{%
\prime}t^{c}x^{d}=A^{2}y^{a}z^{b+1}t^{c-1}t_{-1}x^{d}+2y^{a}x^{d}y^{%
\prime}z^{b}t^{c-1}+A^{-2}y^{a}x^{d}yz^{b}t^{c-1}x.
\end{equation*}

\item[$(7)$] If $w=y^{a}z^{b}z^{\prime }t^{c}t^{\prime }x^{d}$, then in the 
\emph{KBSM} the word $w$ satisfies the following relation: 
\begin{equation*}
y^{a}z^{b}z^{\prime
}t^{c}t^{\prime}x^{d}=A^{2}y^{a}z^{b+1}t^{c+1}x^{d}+2y^{a}x^{d}yz^{b}t^{c}
+A^{-2}y^{a}x^{d}y_{-1}z^{b}t^{c}x.
\end{equation*}
\end{description}

\noindent Thus, in each of the six cases $(2)-(7)$ we can express the
reduced word $w$ in the form $A^{2}P+2Q+A^{-2}R$ for some appropriate
diagrams (different in each case) $P$, $Q$ and $R$. We define the \emph{%
quasi-final} Kauffman bracket by setting 
\begin{equation*}
\left\langle w\right\rangle _{qf}=A^{2}\left\langle \left\langle
P\right\rangle _{rr}\right\rangle _{qf}+2\left\langle \left\langle
Q\right\rangle _{rr}\right\rangle _{qf}+A^{-2}\left\langle \left\langle
R\right\rangle _{rr}\right\rangle _{qf}.
\end{equation*}
To see that this inductive definition terminates resulting with the linear
combination of the quasi-final words, we notice that in $Q$ and $R$ the sum
of components of types $y$, $z$ and $t$ is decreased by one. Moreover, for $%
P $ an arrow is moved from the $y$-type component to the component of $z$-
or $t$-types, and from component of the $z$-type to component of $t$-type,
yielding finally a quasi-final diagram.

\begin{definition}
\label{def_kbqf}Let $D$ be a diagram and let $\left\langle D\right\rangle
_{rr}$ be a linear combination of some reduced words $w$, that is$,$ $%
\left\langle D\right\rangle _{rr}=\sum_{w}P_{w}w$ for $P_{w}\in R$. We
define the \textit{quasi-final} Kauffman bracket of $D$ by setting$:$

\begin{equation*}
\left\langle D\right\rangle_{qf}=\sum_{w}P_{w}\left\langle w\right\rangle_{qf}.
\end{equation*}
\end{definition}

\begin{figure}[h] 
   \centering
   \includegraphics[scale=0.6]{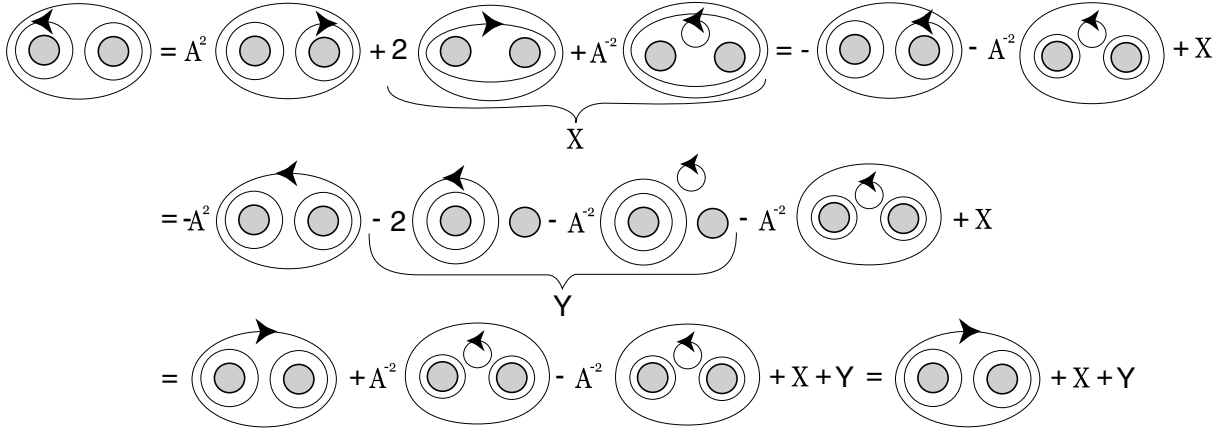} 
   \caption{Passing the arrow from the component of $y$-type to one of $t$-type}
   \label{figf_xt}
\end{figure}

\noindent In the \emph{KBSM} the relation shown in Figure~\ref{figf_xt} holds.

\begin{figure}[h] 
   \centering
   \includegraphics[scale=0.6]{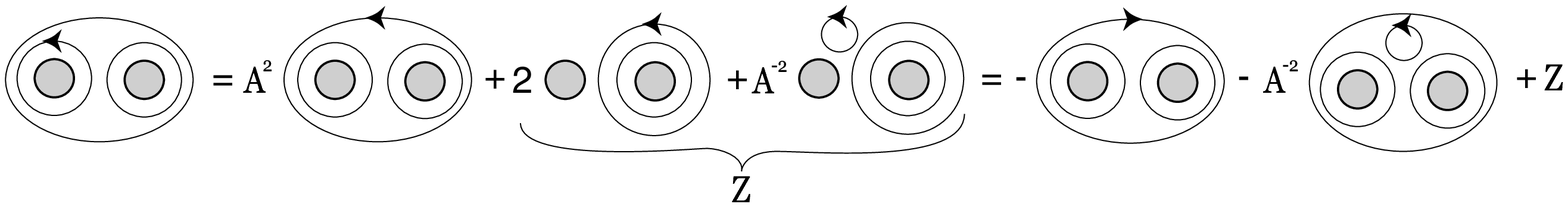} 
   \caption{Passing the arrow from the component of $y$-type to one of $t$-type}
   \label{figf_xt2}
\end{figure}

\noindent By passing the arrow from the component of $y$-type to the $t$%
-type component directly we obtain another relation presented in Figure~\ref{figf_xt2}.

\begin{figure}[h] 
   \centering
   \includegraphics[scale=0.6]{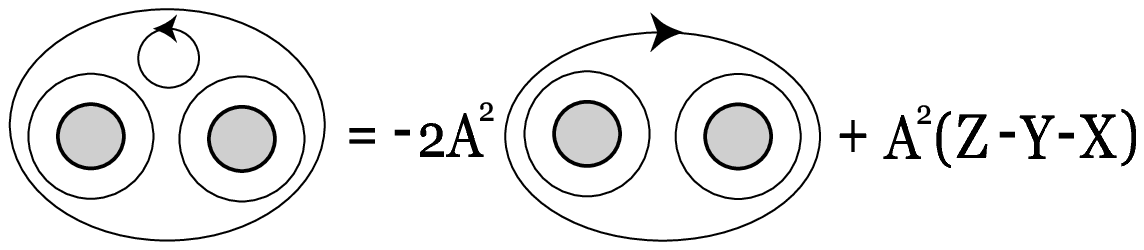} 
   \caption{Eliminating diagrams with all four types of components}
   \label{figf_xt3}
\end{figure}

\noindent The elements $X$, $Y$ and $Z$ have less components of the $y$, $z$
or $t$-types. Setting equal the last terms of the two equations above and
rearranging the terms gives the relation in the \emph{KBSM} as in Figure~\ref{figf_xt3}.

\begin{figure}[h] 
   \centering
   \includegraphics[scale=0.6]{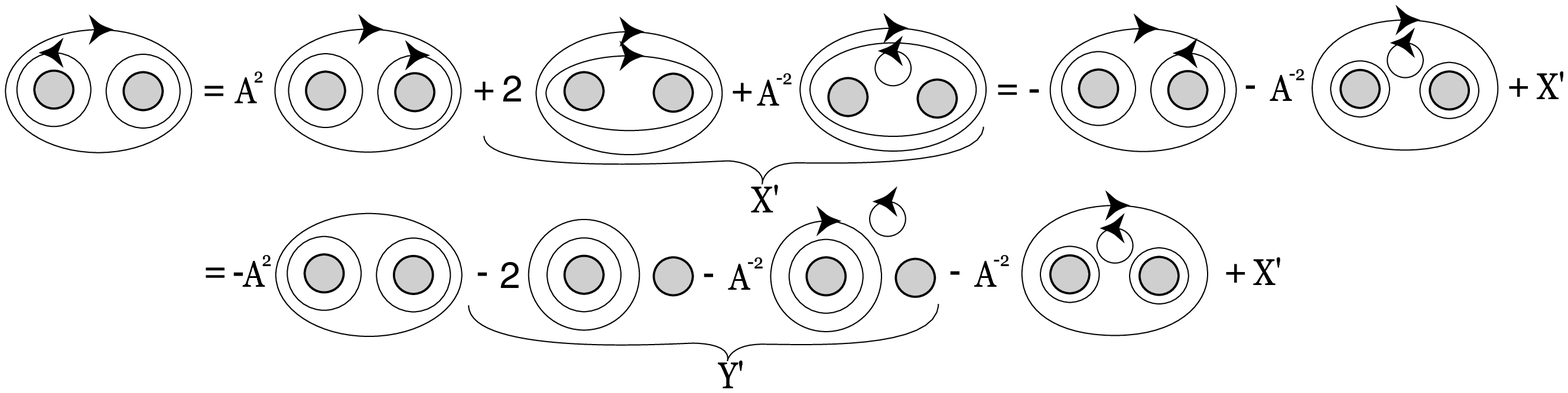} 
   \caption{Passing the arrow from the component of $y$-type to the $t$-type component}
   \label{figf_xtp}
\end{figure}

\begin{figure}[h] 
   \centering
   \includegraphics[scale=0.6]{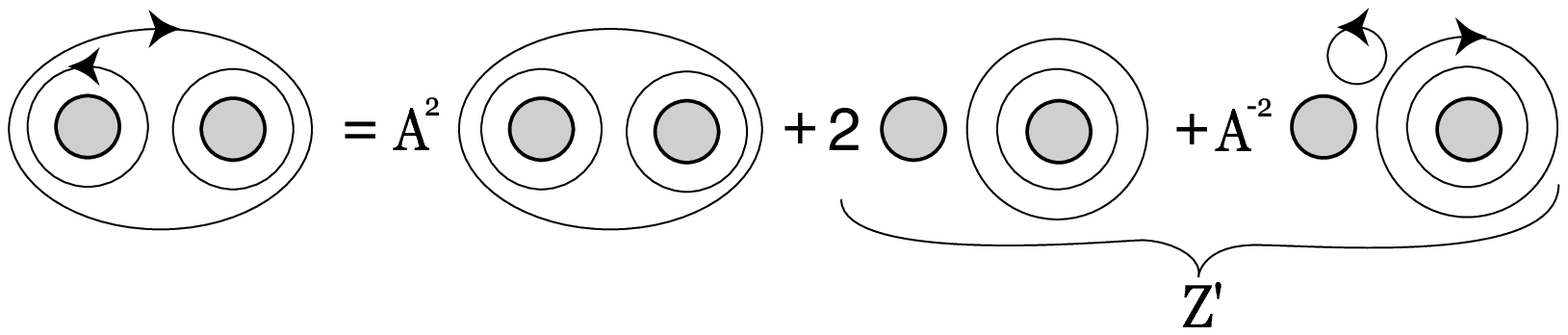} 
   \caption{Passing the arrow from the component of $y$-type to the $t$-type component}
   \label{figf_xtp2}
\end{figure}

\noindent Analogously, when there is an arrow on the $t$-type component,
there are two ways of moving the arrow from the $y$ to the $t$-type
component: first via $z$ as shown in Figure~\ref{figf_xtp} or directly as
shown in Figure~\ref{figf_xtp2}.

\begin{figure}[h] 
   \centering
   \includegraphics[scale=0.6]{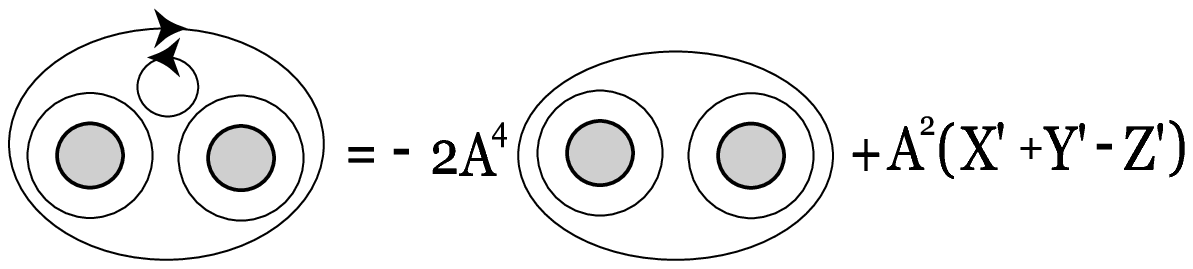} 
   \caption{Eliminating diagrams with all four types of components}
   \label{figf_xtp3}
\end{figure}

\noindent Setting equal the last terms of these two equations and
rearranging terms gives a relation in the \emph{KBSM} that is shown in
Figure~\ref{figf_xtp3}.

\noindent A reduced word is called \emph{final} if it is quasi-final and
does not contain components of all the $4$ types ($x$, $y$, $z,$ and $t$).
We define inductively the final bracket $\left\langle \hspace{0.2cm}%
\right\rangle _{f}$ as follows:\medskip

\noindent $(1)$ If $w$ is final then we set $\left\langle w\right\rangle
_{f}=w$.\medskip

\noindent If $w$ is quasi-final but not final, it has the form of the \emph{%
LHS} shown in Figure~\ref{figf_xt3} or Figure~\ref{figf_xtp3}. Define $%
\left\langle w\right\rangle _{f}$ inductively to be $\left\langle \hspace{%
0.2cm}\right\rangle _{f}$ applied to the \emph{RHS} of the corresponding
equation shown either in Figure~\ref{figf_xt3} or Figure~\ref{figf_xtp3} in
terms of formulas.\medskip

\noindent $(2)$ If $w=y^{a}z^{b}t^{c}x^{d+1}$ then

\begin{equation*}
\left\langle w\right\rangle _{f}=-2A^{2}\left\langle
y^{a}z^{b}t^{c-1}t^{\prime }x^{d}\right\rangle _{f}+A^{2}\left\langle
\left\langle Z-Y-X\right\rangle_{qf}\right\rangle_{f},
\end{equation*}%
where

\begin{equation*}
\begin{tabular}{ll}
$Z=2y^{a-1}z^{b}x^{d}z^{\prime }t^{c-1}+A^{-2}y^{a-1}z^{b}x^{d}zt^{c-1}x$ & $%
Y=-2y^{a}x^{d}y^{\prime}z^{b-1}t^{c-1}-A^{-2}y^{a}x^{d}yz^{b-1}t^{c-1}x$ \\ 
$X=2y^{a-1}z^{b-1}t^{c}x^{d}t^{\prime }+A^{-2}y^{a-1}z^{b-1}t^{c}x^{d}tx.$ & 
\end{tabular}%
\end{equation*}
\noindent $(3)$ If $w=y^{a}z^{b}t^{c}t^{\prime }x^{d+1}$ then 
\begin{equation*}
\left\langle w\right\rangle_{f}=-2A^{4}\left\langle
y^{a}z^{b}t^{c+1}x^{d}\right\rangle_{f}+A^{2}\left\langle \left\langle
Y^{\prime }+X^{\prime }-Z^{\prime }\right\rangle_{qf}\right\rangle_{f}
\end{equation*}
where 
\begin{equation*}
\begin{tabular}{ll}
$Z^{\prime}=2y^{a-1}z^{b}x^{d}zt^{c}+A^{-2}y^{a-1}z^{b}x^{d}z_{-1}t^{c}x$ & $%
Y^{\prime}=-2y^{a}x^{d}y z^{b-1}t^{c}-A^{-2}y^{a}x^{d}y_{-1}z^{b-1}t^{c}x$
\\ 
$X^{\prime }=2y^{a-1}z^{b-1}t^{c}t^{\prime
}x^{d}t^{\prime}+A^{-2}y^{a-1}z^{b-1}t^{c}t^{\prime }x^{d}tx.$ & 
\end{tabular}%
\end{equation*}
\noindent In $X$, $Y$, $Z$, $X^{\prime }$, $Y^{\prime }$ and $Z^{\prime }$
the total number of $y$-, $z$- and $t$-types decreases, and in the other
terms the number of the components of the $x$-type decreases, so the
induction results with words that are final.

\begin{definition}
\label{def_kbff} Let $D$ be a diagram and $\left\langle D\right\rangle _{qf}$
be a linear combination of some quasi-final words $w$, that is, $%
\left\langle D\right\rangle _{qf}=\sum_{w}P_{w}w$ for some $P_{w}\in R$.
Define the \textit{final} Kauffman bracket of $D$ by$:$%
\begin{equation*}
\left\langle D\right\rangle _{f}=\sum_{w}P_{w}\left\langle w\right\rangle
_{f}.
\end{equation*}
\end{definition}

\noindent \textbf{Remark 1\quad }\label{qf_f}By definition, the quasi-final
bracket $\left\langle \hspace*{0.2cm}\right\rangle _{qf}$ preserves all the
equalities appearing in Figure~\ref{figf_xt}. However, $\left\langle 
\hspace*{0.2cm}\right\rangle _{qf}$ does not preserve the first equality in
Figure~\ref{figf_xt2}. But the final bracket $\left\langle \hspace*{0.2cm}%
\right\rangle _{f}$ preserves this equality since it satisfies, by its
definition, the relation shown in Figure~\ref{figf_xt3}. Analogously, $%
\left\langle \hspace*{0.2cm}\right\rangle _{qf}$ preserves all equalities in
Figure~\ref{figf_xtp} and does not preserve the equality in Figure~\ref%
{figf_xtp2}. However, again $\left\langle \hspace*{0.2cm}\right\rangle _{f}$
preserves this equality since, by the definition, it satisfies the relation
in Figure~\ref{figf_xtp3}.\medskip

\noindent To show the invariance of $\left\langle \hspace*{0.2cm}%
\right\rangle _{f}$ under $\Omega _{5}$-moves of all three types, first note
that Lemma~\ref{xinout2} can clearly be extended from $A\times S^{1}$ to the
case of $F_{0,3}\times S^{1}$ (one just considers the components of $y$, $z$
or $t$-type separately).

\begin{theorem}
\label{MainResult} The final Kauffman bracket $\left\langle \hspace*{0.2cm}%
\right\rangle _{f}$ is invariant under all regular Reidemeister moves.
Therefore $S_{2,\infty }(F_{0,3}\times S^{1};R,A)$ is a free $R$-module with
basis that consists of all final words.
\end{theorem}

\noindent 
\begin{proof}%
As before, we assume that, except for the crossing at which $\Omega _{5}$%
-move is to be applied, all other crossings for both diagrams $D$ and $%
D^{\prime }$ were smoothed and all trivial components were already removed.
Let $p$ be the number of components of $y$, $z$ or $t$-types in $D$ (and
also in $D^{\prime }$) without counting the component with the crossing. The
proof of invariance of $\left\langle \hspace*{0.2cm}\right\rangle _{f}$
under the $\Omega _{5}$-move is done by induction on the number $p$. We
first show this invariance for moves of types $\mathbf{I}$, $\mathbf{II}$
and $\mathbf{III}$ when $p=0$.\medskip

\noindent Consider the $\Omega _{5}$-move of type $\mathbf{I}$ shown in
Figure~\ref{figqf_5a0}. Applying Lemma~\ref{xinout2}, we push all $x$-type
components and all arrows out of the component that is involved in the $%
\Omega _{5}$-move, thus it is sufficient to consider the situation when
there are $d$ of $x$-type components outside, where $d\in \mathbb{N}$ and $m$
and $n$ are $0$ or $1$.

\begin{figure}[h] 
   \centering
   \includegraphics[scale=0.5]{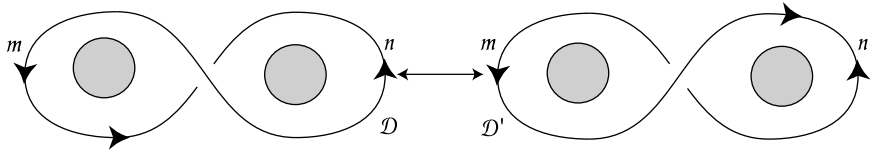} 
   \caption{Invariance of $\left\langle 
\hspace{0.2cm}\right\rangle _{f}$ under $\Omega _{5}$ of type $\mathbf{I}$}
   \label{figqf_5a0}
\end{figure}

\noindent If $m=n=0$ then using relation~\ref{eq1}, which holds already for $%
\left\langle \hspace*{0.2cm}\right\rangle _{srr}$, we have:%
\begin{eqnarray*}
\left\langle D\right\rangle _{f} &=&A\left\langle x^{d}t_{-1}\right\rangle
_{f}+A^{-1}\left\langle y^{\prime }zx^{d}\right\rangle
_{f}=-A^{-1}\left\langle x^{d}t^{\prime }\right\rangle
_{f}-A^{-3}\left\langle x^{d}tx\right\rangle _{f}+A^{-1}\left\langle
y^{\prime }zx^{d}\right\rangle _{f} \\
&=&-A^{-1}\left\langle x^{d}t^{\prime }\right\rangle _{f}-A^{-3}\left\langle
x^{d}tx\right\rangle _{f}+A\left\langle yz_{-1}x^{d}\right\rangle
_{f}+2A^{-1}\left\langle x^{d}t^{\prime }\right\rangle
_{f}+A^{-3}\left\langle x^{d}tx\right\rangle _{f} \\
&=&A^{-1}\left\langle x^{d}t^{\prime }\right\rangle _{f}+A\left\langle
yz_{-1}x^{d}\right\rangle _{f} \\
&=&\left\langle D^{\prime }\right\rangle _{f}.
\end{eqnarray*}%
Analogously, for $m=0$ and $n=1$ we have%
\begin{eqnarray*}
\left\langle D\right\rangle _{f} &=&A\left\langle x^{d}t_{-2}\right\rangle
_{f}+A^{-1}\left\langle y^{\prime }z^{\prime }x^{d}\right\rangle _{f} \\
&=&-A^{-1}\left\langle x^{d}t\right\rangle _{f}-A^{-3}\left\langle
x^{d}t_{-1}x\right\rangle _{f}+A\left\langle yzx^{d}\right\rangle
_{f}+2A^{-1}\left\langle x^{d}t\right\rangle _{f}+A^{-3}\left\langle
x^{d}t_{-1}x\right\rangle _{f} \\
&=&A\left\langle yzx^{d}\right\rangle _{f}+A^{-1}\left\langle
x^{d}t\right\rangle _{f} \\
&=&\left\langle D^{\prime }\right\rangle _{f}
\end{eqnarray*}%
If $m=1$ and $n=0$, we have%
\begin{eqnarray*}
\left\langle D\right\rangle _{f} &=&A\left\langle x^{d}t_{-2}\right\rangle
_{f}+A^{-1}\left\langle y_{2}zx^{d}\right\rangle _{f}=-A^{-1}\left\langle
x^{d}t\right\rangle _{f}-A^{-3}\left\langle x^{d}t_{-1}x\right\rangle
-A\left\langle yzx^{d}\right\rangle _{f}-A^{-3}\left\langle y^{\prime
}zx^{d+1}\right\rangle _{f} \\
&=&-A\left\langle yzx^{d}\right\rangle _{f}-2A^{-1}\left\langle
x^{d}t\right\rangle _{f}-A^{-3}\left\langle x^{d}t_{-1}x\right\rangle
_{f}-A^{-3}\left\langle y^{\prime }zx^{d+1}\right\rangle
_{f}+A^{-1}\left\langle x^{d}t\right\rangle _{f} \\
&=&-A^{-1}\left\langle y^{\prime }z^{\prime }x^{d}\right\rangle
_{f}-A^{-3}\left\langle y^{\prime }zx^{d+1}\right\rangle
_{f}+A^{-1}\left\langle x^{d}t\right\rangle _{f}=A\left\langle y^{\prime
}z_{-1}x^{d}\right\rangle _{f}+A^{-1}\left\langle x^{d}t\right\rangle _{f} \\
&=&\left\langle D^{\prime }\right\rangle _{f}
\end{eqnarray*}%
If $m=1$ and $n=1,$ we have%
\begin{equation*}
\hspace*{-0.1cm}%
\begin{array}{ccl}
\hspace*{-0.1cm}\left\langle D\right\rangle _{f}\hspace*{-0.1cm} & = & 
A\left\langle x^{d}t_{-3}\right\rangle _{f}+A^{-1}\left\langle
y_{2}z^{\prime }x^{d}\right\rangle _{f}=-A^{-1}\left\langle
x^{d}t_{-1}\right\rangle _{f}-A^{-3}\left\langle x^{d}t_{-2}x\right\rangle
_{f}-A\left\langle yz^{\prime }x^{d}\right\rangle _{f}-A^{-3}\left\langle
y^{\prime }z^{\prime }x^{d+1}\right\rangle _{f}\smallskip \\ 
& = & -A^{-1}\left\langle x^{d}t_{-1}\right\rangle _{f}+A^{-5}\left\langle
x^{d}tx\right\rangle _{f}+A^{-7}\left\langle x^{d}t_{-1}x^{2}\right\rangle
_{f}-A\left\langle yz^{\prime }x^{d}\right\rangle _{f}-A^{-1}\left\langle
yzx^{d+1}\right\rangle _{f}\smallskip \\ 
& - & 2A^{-3}\left\langle x^{d+1}t\right\rangle _{f}-A^{-5}\left\langle
x^{d+1}t_{-1}x\right\rangle _{f}=-A^{-1}\left\langle
x^{d}t_{-1}\right\rangle _{f}+A^{-5}\left\langle x^{d}tx\right\rangle
_{f}+A^{-7}\left\langle x^{d}t_{-1}x^{2}\right\rangle _{f}\smallskip \\ 
& + & A^{3}\left\langle yz_{-1}x^{d}\right\rangle _{f}-2A^{-3}\left\langle
x^{d+1}t\right\rangle _{f}+A^{-1}\left\langle x^{d}tx\right\rangle
_{f}-A^{-5}\left\langle x^{d}tx\right\rangle _{f}-A^{-7}\left\langle
x^{d}t_{-1}x^{2}\right\rangle _{f}\smallskip \\ 
& = & -A^{-1}\left\langle x^{d}t_{-1}\right\rangle _{f}+A^{3}\left\langle
yz_{-1}x^{d}\right\rangle _{f}+2A\left\langle x^{d}t^{\prime }\right\rangle
_{f}-2A^{-3}\left\langle x^{d}t^{\prime }\right\rangle
_{f}-2A^{-5}\left\langle x^{d}tx\right\rangle _{f}+A^{-1}\left\langle
x^{d}tx\right\rangle _{f}\smallskip \\ 
& = & A^{3}\left\langle yz_{-1}x^{d}\right\rangle _{f}+2A\left\langle
x^{d}t^{\prime }\right\rangle _{f}+A^{-1}\left\langle x^{d}tx\right\rangle
_{f}-A^{-1}\left\langle x^{d}t_{-1}\right\rangle _{f}+2A^{-1}\left\langle
x^{d}t_{-1}\right\rangle _{f}\smallskip \\ 
& = & A\left\langle y^{\prime }zx^{d}\right\rangle _{f}+A^{-1}\left\langle
x^{d}t_{-1}\right\rangle _{f}=\left\langle D^{\prime }\right\rangle _{f}.%
\end{array}%
\end{equation*}%
Now, if $p=0$, the situation is similar for the $\Omega _{5}$-moves of types 
$\mathbf{II}$ and $\mathbf{III}$. In the formulas one just has to permute $z$
with $t$, keeping $y$ (for the type $\mathbf{III)}$, or permute $y$ with $z$
and $z$ with $t$ (for the type $\mathbf{II}$). Note also that in these cases
quasi-final bracket $\left\langle \hspace*{0.2cm}\right\rangle _{qf}$ is
unchanged just like the final bracket $\left\langle \hspace*{0.2cm}%
\right\rangle _{f}$.

\noindent Now, by induction, let us assume that the final bracket $%
\left\langle \hspace*{0.2cm}\right\rangle _{f}$ is invariant under the $%
\Omega _{5}$-moves of types $\mathbf{I}$, $\mathbf{II},$ and $\mathbf{III}$
that involve less than $p$ components of $y$, $z$ or $t$-type (counting
without the component with the crossing). Let $D$ and $D^{\prime }$ have $p$
such components. The situation for the $\Omega _{5}$-move of type $\mathbf{I}
$ is shown in Figure~\ref{figqf_5a}.

\begin{figure}[h] 
   \centering
   \includegraphics[scale=0.45]{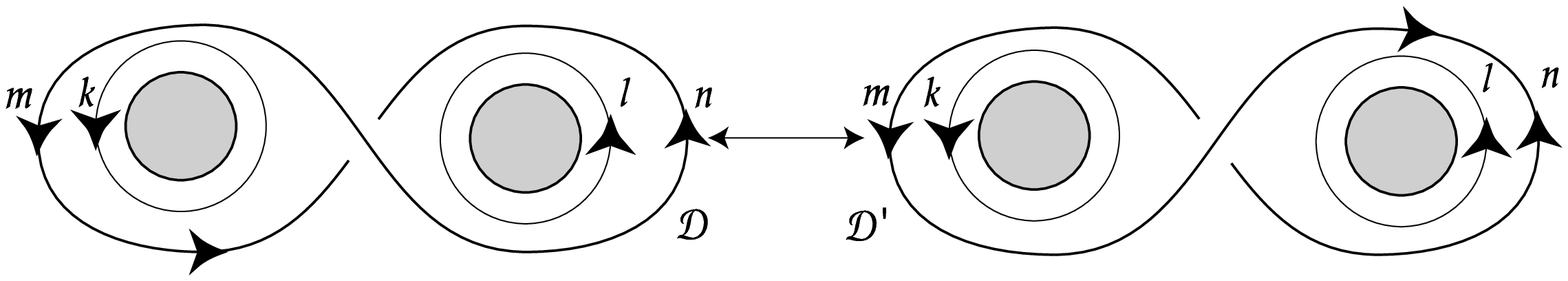} 
   \caption{Invariance of $\left\langle \hspace*{0.2cm}\right\rangle _{f}$ under the $\Omega _{5}$-move of type $\mathbf{I}$}
   \label{figqf_5a}
\end{figure}

\noindent We may assume that the bracket $\left\langle \hspace*{0.2cm}%
\right\rangle _{rr}$ is applied until all possible arrows appear only on the
most external component of the type $t$ (not shown in Figure~\ref{figqf_5a})
and on the components shown in Figure~\ref{figqf_5a}. Using Lemma~\ref%
{xinout2} we also assume that $k$ and $l$ are equal to $0$ or $1$. If both
of them are equal to $0$ then we can use the same argument as in the case $%
p=0$ since the interior components of the $y$-type and $z$-type and all the
components of the $t$-type do not play any role in the previous calculation.
Suppose now that $k=1$ and $l=0$. In that case we apply the $\Omega _{2}$%
-move followed by the $\Omega _{5}$-move to both diagrams $D$ and $D^{\prime
}$ to obtain a situation when we can show easily that $\left\langle \hspace*{%
0.2cm}\right\rangle _{f}$ is the same for both diagrams. These moves are
shown in Figure~\ref{figqf_5aa}.

\begin{figure}[h] 
   \centering
   \includegraphics[scale=0.6]{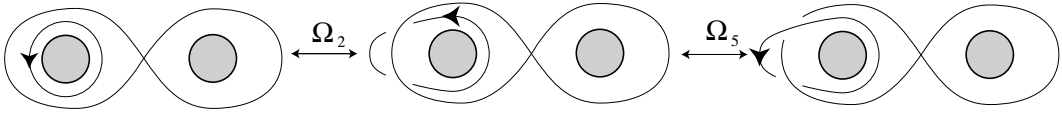} 
   \caption{Changing $D$ and $D^{\prime }$ by $\Omega _{2}$ and $\Omega _{5}$-moves}
   \label{figqf_5aa}
\end{figure}

\noindent The $\Omega _{2}$-move does not change any of the brackets. The
effect of applying the $\Omega _{5}$-move can be analyzed by looking at all
possible smoothings of crossings except the one at which the $\Omega _{5}$%
-move is applied. For all such smoothings we obtain three moves under which $%
\left\langle \hspace*{0.2cm}\right\rangle _{rr}$ is invariant and one move
of type $\mathbf{III}$ but with one less component of $y$, $z$ or $t$-type,
so by the induction hypothesis $\left\langle \hspace*{0.2cm}\right\rangle
_{f}$ is invariant under this move of type $\mathbf{III}$. Thus, if $D_{1}$
is obtained from $D$ by the application of the two Reidemeister moves
mentioned before and, in a similar way, $D_{1}^{\prime }$ is obtained from $%
D^{\prime }$, we have%
\begin{equation*}
\left\langle D\right\rangle _{f}=\left\langle D_{1}\right\rangle _{f}\text{
and }\left\langle D^{\prime }\right\rangle _{f}=\left\langle D_{1}^{\prime
}\right\rangle _{f}.
\end{equation*}%
Now the $\Omega _{5}$-move between $D_{1}$ and $D_{1}^{\prime }$ is again
expressed by smoothing the two crossings not involved in the $\Omega _{5}$%
-move yielding three $\Omega _{5}$-moves under which $\left\langle \hspace*{%
0.2cm}\right\rangle _{rr}$ is invariant and the $\Omega _{5}$-move of type $%
\mathbf{I}$ for which $k=0$ instead of $1$ (and $l$ remains equal to $0$).
Thus we have $\left\langle D_{1}\right\rangle _{f}=\left\langle
D_{1}^{\prime }\right\rangle _{f}$ using the proof of the preceding case $%
k=l=0$; therefore, it follows that $\left\langle D\right\rangle
_{f}=\left\langle D^{\prime }\right\rangle _{f}$. The same proof works for $%
k=0$ and $l=1$. Having established these cases, one applies analogous
arguments for the case $k=l=1$.\medskip

\noindent Consider now the $\Omega _{5}$-move of type $\mathbf{II}$ between
diagrams $D$ and $D^{\prime }$ which involve $p$ components of $y$, $z$ or $%
t $-type. The situation is shown in Figure~\ref{figqf_5b}.

\begin{figure}[h] 
   \centering
   \includegraphics[scale=0.8]{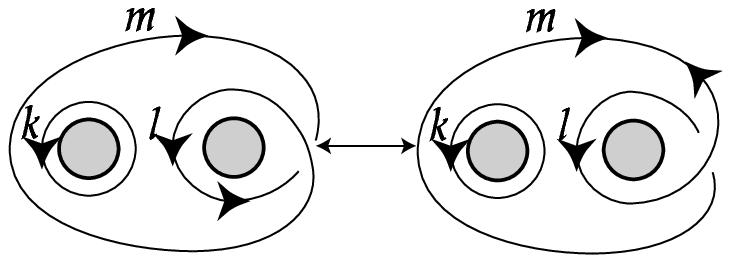} 
   \caption{Invariance of $\left\langle \hspace*{0.2cm}\right\rangle _{f}$ under the $\Omega _{5}$-move of type $\mathbf{II}$}
   \label{figqf_5b}
\end{figure}

\noindent Applying $\left\langle \hspace*{0.2cm}\right\rangle _{rr}$ one
arrives at the case $k=0$ or $k=1$. If $k=0$, then the proof is as before:
if there are no extra arrows except on the $x$-type components and the ones
shown in Figure~\ref{figqf_5b}, then the situation is as in the case $p=0$;
if there are such arrows, then the proof is like for the case of $\Omega
_{5} $-move of type $\mathbf{I}$, where one used two Reidemeister moves to
decompose the $\Omega _{5}$-move into moves for which $\left\langle \hspace*{%
0.2cm}\right\rangle _{f}$ is invariant.\medskip

\noindent If $k=1$ the similar arguments can be applied again; the situation
is shown in Figure~\ref{figqf_5bb}.

\begin{figure}[h] 
   \centering
   \includegraphics[scale=0.75]{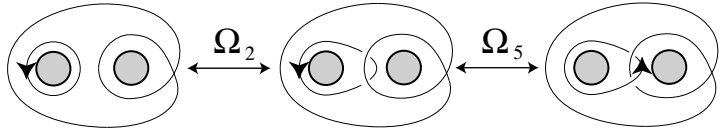} 
   \caption{Changing $D$ and $D^{\prime }$ by $\Omega _{2}$ and $\Omega _{5}$ moves}
   \label{figqf_5bb}
\end{figure}

\noindent Namely, we observe that the $\Omega _{2}$-move leaves the bracket $%
\left\langle \hspace*{0.2cm}\right\rangle _{f}$ unchanged for both $D$ and $%
D^{\prime }$. The $\Omega _{5}$-move is expressed by smoothing the crossings
yielding three moves for which $\left\langle \hspace*{0.2cm}\right\rangle
_{rr}$ is preserved and one $\Omega _{5}$-move of type $\mathbf{I}$ with the
number $p$ unchanged, for which we already showed that the bracket $%
\left\langle \hspace*{0.2cm}\right\rangle _{f}$ is unchanged. Finally, after
these two Reidemeister moves, the original $\Omega _{5}$-move of type $%
\mathbf{II}$ is expressed using three moves for which the bracket $%
\left\langle \hspace*{0.2cm}\right\rangle _{rr}$ is unchanged and the $%
\Omega _{5}$-move of type $\mathbf{II}$ for which $k=0$.

\noindent It remains to show the invariance of $\left\langle \hspace*{0.2cm}%
\right\rangle _{f}$ under the $\Omega _{5}$-move of type $\mathbf{III}$;
such a case is illustrated in Figure~\ref{figqf_5c}.

\begin{figure}[h] 
   \centering
   \includegraphics[scale=0.75]{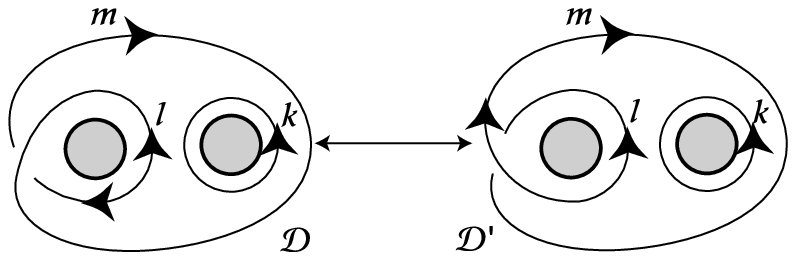} 
   \caption{Invariance of $\left\langle \hspace*{0.2cm}\right\rangle _{f}$ under $\Omega _{5}$ move of type $\mathbf{III}$}
   \label{figqf_5c}
\end{figure}

\noindent Again, there are few cases to consider. First, if there are no
components of the type $z$ (i.e. neither $z$ nor $z^{\prime }$ present in
the words), then the proof is as before: either there are no arrows except
on the $x$-type components and the ones shown in Figure~\ref{figqf_5c}, in
which case the calculation is the same as for the case $p=0$; or there are
such arrows on the components of $y$ or $t$-type, in which case the proof is
the same as it was in that case for $\Omega _{5}$-move of types $\mathbf{I}$
and $\mathbf{II}$.\medskip

\noindent It remains to check the invariance of $\left\langle \hspace*{0.2cm}%
\right\rangle _{f}$ when the $z$-type components appear in both diagrams $D$
and $D^{\prime }$. Applying Lemma~\ref{xinout2}, we can assume that $k$, $l,$
and $m$ are equal to $0$ or $1$, and all $x$-type components are in between
the components of the three other types and there are no extra
arrows.\medskip

\noindent Suppose that $k=0$. In the computations of $\left\langle
D\right\rangle _{f}$ and $\left\langle D^{\prime }\right\rangle _{f}$ one
proceeds as in the case $p=0$ except for the situations when an arrow
appearing on a $y$-type component is moved to become an arrow on a $t$-type
component. By definition of $\left\langle \hspace*{0.2cm}\right\rangle _{qf}$
this has to be done via a $z$-type component. However, by Remark 1, for $%
\left\langle \hspace*{0.2cm}\right\rangle _{f}$ the same result is obtained
if the arrow is moved directly from the $y$-type component to the $t$-type
component. Thus, for the final bracket $\left\langle \hspace*{0.2cm}%
\right\rangle _{f}$ the component of the $z$-type plays no role in the
calculations which are then like for the case $p=0$.\medskip

\noindent The final case to check is when $k=1$. Again, we decompose the $%
\Omega _{5}$-move of type $\mathbf{III}$ into other moves for which the
invariance of $\left\langle \hspace*{0.2cm}\right\rangle _{f}$ has already
been established. This is shown in Figure~\ref{figqf_5cc}.

\begin{figure}[h] 
   \centering
   \includegraphics[scale=0.7]{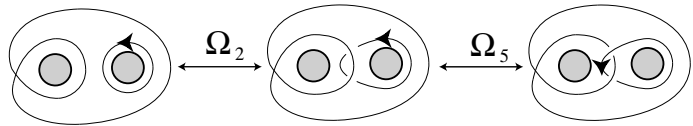} 
   \caption{Changing $D$ and $D^{\prime }$ by $\Omega _{2}$ and $\Omega _{5}$ moves}
   \label{figqf_5cc}
\end{figure}

\noindent The $\Omega _{5}$-move shown in Figure~\ref{figqf_5cc} does not
change $\left\langle \hspace*{0.2cm}\right\rangle _{f}$ since it can be
decomposed by smoothing the crossings not involved in the move into three
moves for which $\left\langle \hspace*{0.2cm}\right\rangle _{rr}$ is
invariant and a $\Omega _{5}$-move of type $\mathbf{I}$. After this, the
desired $\Omega _{5}$-move is applied and, considering the smoothings of the
crossings that are not involved in the move, it can be expressed using three
moves (for which $\left\langle \hspace*{0.2cm}\right\rangle _{rr}$ is
invariant) and the $\Omega _{5}$-move of type $\mathbf{III}$ with $k=0$ and
for which, as it was shown above, the bracket $\left\langle \hspace*{0.2cm}%
\right\rangle _{f}$ is invariant. This finishes our proof.%
\end{proof}%
\medskip

\noindent \textbf{Acknowledgments}

\noindent The authors would like to thank Marie Reed for helpful editorial suggestions.

\end{document}